\documentclass[10pt]{article}
\usepackage[T1]{fontenc}
\usepackage{lmodern}
\usepackage[margin=1in]{geometry}
\usepackage{amsmath,amssymb,amsthm,mathtools}
\usepackage{enumitem,microtype,etoolbox,needspace}
\usepackage{indentfirst}
\usepackage{tikz}
\usetikzlibrary{arrows.meta}
\tikzset{
  rt picture/.style={x=1cm,y=1cm,font=\small,
    line cap=round,line join=round},
  rt vertex/.style={circle,draw=black,fill=white,line width=.65pt,
    minimum size=4pt,inner sep=0pt,outer sep=0pt},
  rt actual/.style={rt vertex,fill=black},
  rt formal/.style={draw=blue!55!black,line width=.9pt,
    dash pattern=on 3pt off 2pt},
  rt represent/.style={draw=black,line width=1pt},
  rt portion/.style={draw=black!45,line width=.85pt},
  rt map/.style={-{Latex[length=2mm,width=1.4mm]},
    draw=black!60,line width=.7pt},
  rt shade/.style={rounded corners=3pt,fill=black!4,draw=none},
  rt note/.style={font=\footnotesize,align=left}
}
\usepackage[font=small]{caption}

\usepackage[colorlinks=true,linkcolor=blue,citecolor=blue,urlcolor=blue]{hyperref}
\hypersetup{pdftitle={Chv\'atal's Toughness Conjecture for Dense Almost Regular Graphs},pdfauthor={Songling Shan}}
\setlist[enumerate]{leftmargin=*,itemsep=3pt,topsep=5pt}
\newtheorem{theorem}{Theorem}[section]
\newtheorem{lemma}[theorem]{Lemma}
\newtheorem{proposition}[theorem]{Proposition}
\newtheorem{conjecture}[theorem]{Conjecture}
\newtheorem{claim}{Claim}
\newcounter{claimblock}

\newenvironment{claimproof}[1]
  {\begin{proof}[Proof of Claim~\ref{#1}]}
  {\end{proof}}
\BeforeBeginEnvironment{theorem}{\Needspace{6\baselineskip}}
\BeforeBeginEnvironment{lemma}{\Needspace{6\baselineskip}}
\BeforeBeginEnvironment{claim}{\Needspace{6\baselineskip}}
\BeforeBeginEnvironment{proposition}{\Needspace{6\baselineskip}}
\newcommand{\RN}{\operatorname{RN}}
\newcommand{\cD}{\mathcal D}
\newcommand{\cU}{\mathcal U}
\newcommand{\cP}{\mathcal P}
\newcommand{\cV}{\mathcal V}
\newcommand{\cH}{\mathcal H}
\newcommand{\cT}{\mathcal T}
\newcommand{\bP}{\mathbf P}
\newcommand{\ol}[1]{\overline{#1}}
\title{Chv\'atal's Toughness Conjecture for Dense Almost Regular Graphs}
\author{Songling Shan\thanks{Auburn University, Department of Mathematics and Statistics,
Auburn, AL 36849, USA. Email: \texttt{szs0398@auburn.edu}. Supported in part
by NSF grant DMS-2451895.}}
\date{}

\begin{document}
\maketitle

\begin{abstract}
Chv\'atal's toughness conjecture asserts that there is a constant
$t_0$ such that every $t_0$-tough graph on at least three vertices is
Hamiltonian. Most partial results on this conjecture concern graphs
with a restricted structure, such as chordal graphs and $H$-free graphs
for small linear forests $H$. We confirm the conjecture with $t_0=4$
for dense almost regular graphs, without any structural assumption:
for every $\alpha>0$, there exist $\varepsilon>0$ and $n_0$ such that
every $4$-tough graph $G$ on $n\ge n_0$ vertices with
$\delta(G)\ge\alpha n$ and $\Delta(G)-\delta(G)\le\varepsilon n$ is
Hamiltonian. This result also bears on a conjecture of Bollob\'as and,
independently, H\"aggkvist from the 1970s that every $k$-connected
$D$-regular graph on $n$ vertices with $D\ge n/(k+1)$ is Hamiltonian.
This conjecture is false for every $k\ge4$, even for large $n$. Our
result shows that its
conclusion holds for all large $n$ when $k$-connectivity is replaced by
$4$-toughness, and that the degree condition can then be weakened to
$D\ge\alpha n$ for an arbitrarily small fixed $\alpha>0$. We also
construct $3/2$-tough regular non-Hamiltonian graphs of linear degree,
showing that the constant $4$ cannot be replaced by $3/2$.
\end{abstract}


\section{Introduction}\label{sec:intro}

Let
$t\ge0$ be a real number. A graph $G$ is \emph{$t$-tough} if
$|S|\ge t\,c(G-S)$ for every $S\subseteq V(G)$ with $c(G-S)\ge2$, where
$c(H)$ denotes the number of components of a graph $H$. The concept of
toughness was introduced by Chv\'atal \cite{Chvatal1973} in 1973.
Clearly, every Hamiltonian graph is $1$-tough, but the converse is not
true. Nevertheless, Chv\'atal \cite{Chvatal1973} proposed the following
conjecture.

\begin{conjecture}[Chv\'atal~\cite{Chvatal1973}]\label{conj:chvatal}
There exists a constant $t_0$ such that every $t_0$-tough graph on at
least three vertices is Hamiltonian.
\end{conjecture}

Conjecture~\ref{conj:chvatal} is still open, and Bauer, Broersma, and
Veldman \cite{BBV2000} showed that such a constant $t_0$, if it exists,
is at least $9/4$. Most of the known partial results concern graphs
with a restricted structure. Since every $t$-tough graph that is not
complete is $2t$-connected, the conjecture holds for planar graphs by
Tutte's theorem that $4$-connected planar graphs are
Hamiltonian~\cite{Tutte1956}, and for claw-free graphs by the result of
Kaiser and Vr\'ana~\cite{KaiserVrana2012} that $5$-connected
claw-free graphs of minimum degree at least $6$ are Hamiltonian. It
also holds for chordal
graphs~\cite{CJKL1998,KabelaKaiser2017} and some of their subclasses,
such as planar chordal graphs~\cite{BHT1999}, split
graphs~\cite{KLM1996}, and interval graphs (implicitly
in~\cite{Keil1985}), and for
cocomparability graphs~\cite{DKS1997}. Another line of work concerns
$H$-free graphs, that is, graphs with no induced subgraph isomorphic to
$H$, where $H$ is a small linear forest. Jung~\cite{Jung1978} showed in 1978 that every
$1$-tough $P_4$-free graph on at least three vertices is Hamiltonian,
and much of the recent work has concerned other small linear forests;
see, for example, \cite{BPP2014,OtaSanka2022,Shan2020} for $H=2K_2$,
\cite{GaoShan2022} for $H=P_3\cup2P_1$, \cite{Shan2021} for
$H=P_2\cup P_3$, \cite{Shan2026} for $H=P_4\cup P_1$, and
\cite{OtaSanka2024,ShiShan2022,XuLiZhou2024} for $H=P_2\cup kP_1$, where a
connectivity condition is also imposed. We refer to the survey of
Bauer, Broersma, and Schmeichel~\cite{BBS2006} for an overview of
earlier results.

Degree conditions provide another approach to the conjecture. For
$t>0$, Bauer, Broersma, van den Heuvel, and
Veldman~\cite[Corollary~12]{Bauer1995} proved that a $t$-tough graph $G$
on $n\ge3$ vertices is Hamiltonian if $\delta(G)>n/(t+1)-1$. Sanka and
Shan~\cite[Theorem~1(a)]{SankaShan2024} proved that the same conclusion
holds under the Ore-type condition $\sigma_2(G)>2n/(t+1)-2$, where
$\sigma_2(G)$ is the minimum degree sum of two nonadjacent vertices of
$G$, with $\sigma_2(G)=+\infty$ when $G$ is complete. For $4$-tough
graphs, the first of these results requires $\delta(G)>n/5-1$.

For regular graphs, Hamiltonicity has also been studied under
connectivity conditions. Bollob\'as~\cite[p.~167]{Bollobas1978} and
H\"aggkvist (see~\cite{Jackson1980}) independently proposed the
following conjecture.

\begin{conjecture}[Bollob\'as; H\"aggkvist]\label{conj:BH}
Let $k$ be a positive integer. Every $k$-connected $D$-regular graph on
$n\ge3$ vertices with $D\ge n/(k+1)$ is Hamiltonian.
\end{conjecture}

We state the conjecture as in~\cite{KLOS2016}, adding the condition
$n\ge3$, which only excludes $K_2$ for $k=1$; a footnote
in~\cite{KLOS2016} notes that the formulation of Bollob\'as has the
weaker degree condition $D\ge n/(k+1)-1$. Jackson~\cite{Jackson1980} proved
Conjecture~\ref{conj:BH} for $k=2$. For $k=3$, the conjecture fails for
small $n$: the Petersen graph has $10$ vertices and is $3$-connected and
$3$-regular with $3\ge10/4$, but it is not Hamiltonian; the case $k=3$
is stated in~\cite[Conjecture~1.1]{KLOS2015} for $n\ge13$. K\"uhn, Lo, Osthus, and
Staden proved this case first for $D\ge(1/4+\varepsilon)n$ and large
$n$~\cite{KLOS2015}, and then for $D\ge n/4$ and large
$n$~\cite{KLOS2016}. The
conjecture fails for every $k\ge4$, even for large $n$: a construction
of Jung and, independently, of Jackson, Li, and Zhu (see
\cite[Section~1.2]{KLOS2015}) gives, for every $m$ divisible by four,
an $m/2$-connected $m$-regular graph $G$ on $4m+1$ vertices with a set
$A$ of vertices such that $G-A$ has $|A|+1$ components. Taking $k\ge4$
and $m\ge2k$ divisible by four, the graph $G$ is $k$-connected and $m\ge(4m+1)/(k+1)$.
However, since $|A|/c(G-A)<1$, the graph $G$ is not $1$-tough, and so
$G$ is not Hamiltonian. Moreover, K\"uhn, Lo,
Osthus, and Staden \cite[Theorem~1.4 and Proposition~8.1]{KLOS2015}
showed that, for fixed integers $k\ge1$ and $r\ge2$ and fixed $\xi>0$, every
$k$-connected $D$-regular graph on $n$ vertices with
$D\ge(1/r+\xi)n$ has a cycle of length at least
$\min\{k/(r-1),1-\xi\}n$ when $n$ is large, and that this bound is
asymptotically best possible in general. Thus no fixed connectivity condition
forces a Hamiltonian cycle in $D$-regular graphs when $D$ is a small
linear function of $n$, and the construction above shows that even
connectivity linear in $n$ does not suffice when $D$ is about $n/4$.

In this paper, we show that toughness, unlike connectivity, forces a
Hamiltonian cycle in dense almost regular graphs of every linear
degree. Our main result is the following.

\begin{theorem}\label{thm:main}
For every $\alpha>0$, there exist $\varepsilon>0$ and an integer $n_0$
such that the following holds. Let $G$ be a $4$-tough graph on
$n\ge n_0$ vertices. If
\begin{equation}\label{eq:main}
\delta(G)\ge\alpha n
\qquad\text{and}\qquad
\Delta(G)-\delta(G)\le\varepsilon n,
\end{equation}
then $G$ is Hamiltonian.
\end{theorem}

Theorem~\ref{thm:main} contributes to both conjectures above. First, it
confirms Conjecture~\ref{conj:chvatal} with $t_0=4$ for dense almost
regular graphs. Unlike the results for the graph classes listed above,
it imposes no condition on the structure of $G$; its hypotheses
concern only the degrees. Moreover, for every $\alpha>0$, the condition
$\delta(G)>n/5-1$ in the result of Bauer, Broersma, van den Heuvel, and
Veldman can be replaced by $\delta(G)\ge\alpha n$ for $4$-tough graphs
$G$ with $\Delta(G)-\delta(G)\le\varepsilon n$ and $n\ge n_0$, where
$\varepsilon$ and $n_0$ are given by Theorem~\ref{thm:main}.

Second, Theorem~\ref{thm:main} shows that toughness rescues
Conjecture~\ref{conj:BH}. Applying it with $\alpha=1/(k+1)$, we obtain
that, for every positive integer $k$, every $4$-tough $D$-regular graph
on $n$ vertices with $D\ge n/(k+1)$ is Hamiltonian if $n$ is
sufficiently large. Thus the conclusion of Conjecture~\ref{conj:BH},
which fails for every $k\ge4$, holds for every $k$ and all sufficiently
large $n$ once $k$-connectivity is replaced by $4$-toughness, and the degree condition can even be
weakened to $D\ge\alpha n$ for any fixed $\alpha>0$. Recall that the
counterexamples above have connectivity linear in $n$ but are not
$1$-tough.

Chv\'atal~\cite[Section~5]{Chvatal1973} already constructed an infinite
family of non-Hamiltonian $3/2$-tough graphs. He showed that the
inflation of a non-Hamiltonian $3$-edge-connected cubic graph, obtained
by replacing every vertex by a triangle, is again a non-Hamiltonian
$3$-edge-connected cubic graph and is
$3/2$-tough~\cite[Corollary~5.3]{Chvatal1973}; repeated inflation,
starting from the Petersen graph, gives his examples. Our
construction in Section~\ref{sec:example} replaces the vertices of the
Petersen graph by almost complete graphs instead. It gives such graphs
that are, in addition, dense, and so it shows that toughness $3/2$ does
not suffice in Theorem~\ref{thm:main}, even for regular graphs.

\begin{proposition}\label{prop:petersen}
For every $0<\alpha<1/10$, there are regular non-Hamiltonian $3/2$-tough graphs of arbitrarily large order $n$ with
minimum degree at least $\alpha n$.
\end{proposition}

\medskip
\noindent\textbf{Organization.} In Section~\ref{sec:prelim}, we list 
the tools for the proof. Section~\ref{sec:robust} states the robust
partition theorem and the linkage results for robust expanders that we
use, and proves Lemma~\ref{lem:unbalanced}, which gives spanning
linkages in bipartite robust expanders whose parts have different
sizes.
Section~\ref{sec:matching} proves Lemma~\ref{lem:twotrees}, which
provides an Euler circuit through all the sets of a partition whose
edges form a matching. In
Section~\ref{sec:proof}, we combine these results to prove
Theorem~\ref{thm:main}. In Section~\ref{sec:example}, we prove
Proposition~\ref{prop:petersen}.

\medskip
\noindent\textbf{Conventions and notation.} All graphs considered in
this paper are finite and simple, except for some auxiliary
multigraphs and one auxiliary digraph that will be specified. Let $G$ be a graph. We denote by
$V(G)$ and $E(G)$ the vertex set and edge set of $G$, respectively,
and we let $e(G)=|E(G)|$; for a digraph, $E(G)$ is its set of arcs,
and so $e(G)$ is its number of arcs. For
$v\in V(G)$, we denote by $N_G(v)$ the set of neighbors of $v$ in $G$,
and $d_G(v)$ is the number of edges incident with $v$; $\delta(G)$ and
$\Delta(G)$ are the minimum and maximum degree of $G$. For
$S\subseteq V(G)$, $G[S]$ is the subgraph of $G$
induced by $S$, and $G-S=G[V(G)\setminus S]$; for $F\subseteq E(G)$, the
graph $G-F$ is obtained from $G$ by deleting the edges in $F$, and we
write $G-v$ for $G-\{v\}$. For $A,B\subseteq V(G)$, not necessarily
disjoint, we denote by $E_G(A,B)$ the set of edges of $G$ with one
endvertex in $A$ and the other in $B$, where every edge with both
endvertices in $A\cap B$ is counted only once, and we let
$e_G(A,B)=|E_G(A,B)|$. For $v\in V(G)$ and $S\subseteq V(G)$, we let
$d_G(v,S)=e_G(\{v\},S)$. If $A$ and $B$ are disjoint, we denote by
$G[A,B]$ the bipartite subgraph of $G$ with vertex set $A\cup B$ and
edge set $E_G(A,B)$. A \emph{multigraph} may have parallel edges but
no loops. For a graph or multigraph $H$ and a set $F\subseteq E(H)$,
we denote by $H[F]$ the subgraph of $H$ with edge set $F$ whose
vertex set consists of the endvertices of the edges in $F$.
For integers $p$ and $q$, let
$[p,q]=\{i\in\mathbb Z: p\le i\le q\}$; in particular, we use $[q]$ for $[1,q]$. We omit the subscript $G$ when
the graph is clear from the context.

For constants, the notation $a\ll b$ means that $a$ is chosen
sufficiently small as a function of $b$ and the constants to its right.
More precisely, $0<a\ll b\ll c$ means that there are nondecreasing
functions $g_1$ and $g_2$ from $(0,1)$ to $(0,1)$, where $g_2$ may
depend on $c$, such that the statement holds whenever $b\le g_1(c)$ and
$a\le g_2(b)$, and
similarly for longer hierarchies. The notation $1/m\ll a$ means that
$m$ is sufficiently large after the constants have been fixed, in the
same sense. All constants in these hierarchies are positive. Whenever
a hierarchy is used in a counting or greedy-choice argument, we give
the explicit inequalities that we need.

\section{Preliminaries}\label{sec:prelim}

The proof of Theorem~\ref{thm:main} needs two kinds of tools. We will
first find a partition $\cV=\{V_1,\ldots,V_k,W_1,\ldots,W_\ell\}$ of
$V(G)$, whose members we call \emph{clusters}, such that $G[V]$ has
certain expansion properties for each $V\in\cV$. The first tool concerns linkages inside each cluster of $\cV$
and is provided in Section~\ref{sec:robust}. The second tool concerns
linkages between different clusters of $\cV$ and is provided in
Section~\ref{sec:matching}.

\subsection{Robust expanders and spanning linkages}\label{sec:robust}

The property that we need inside each cluster is the following. Let $p$
be a positive integer. A graph $G$ is \emph{Hamilton
$p$-linked} if for every $2p$ distinct vertices
$x_1,y_1,\ldots,x_p,y_p$ of $G$, there are vertex-disjoint paths
$Q_1,\ldots,Q_p$ in $G$ such that $Q_i$ is an $(x_i,y_i)$-path for each
$i\in[p]$ and $\bigcup_{i=1}^pV(Q_i)=V(G)$. We call such a family of
paths a \emph{spanning linkage} of $G$ for the pairs $\{x_i,y_i\}$,
$i\in[p]$. Thus a graph is Hamilton $1$-linked if and only if it is
Hamiltonian-connected.

Let $0<\nu,\tau<1$, let $G$ be a graph on $n$ vertices, and let
$S\subseteq V(G)$. The \emph{$\nu$-robust neighborhood} of $S$ in $G$ is
\[
\RN_{\nu,G}(S)=\{v\in V(G):d_G(v,S)\ge\nu n\}.
\]
The graph $G$ is a \emph{robust $(\nu,\tau)$-expander} if
$|\RN_{\nu,G}(S)|\ge|S|+\nu n$ for every $S\subseteq V(G)$ with
$\tau n\le|S|\le(1-\tau)n$. If $(A,B)$ is a bipartition of $V(G)$, we
say that $G$ is a \emph{bipartite robust $(\nu,\tau)$-expander with
bipartition $(A,B)$} if
\[
|\RN_{\nu,G}(S)\cap B|\ge |S|+\nu n
\quad\text{for every}\quad
S\subseteq A\ \text{with}\ \tau|A|\le|S|\le(1-\tau)|A|.
\]
The order of $A$ and $B$ is part of this definition. The graph $G$
need not be bipartite here; deleting its edges within $A$ and within
$B$ does not change the condition above. Note also that, for any
$0<\nu'\le\nu\le\tau\le\tau'<1$, every robust $(\nu,\tau)$-expander is
a robust $(\nu',\tau')$-expander, and every bipartite robust
$(\nu,\tau)$-expander with bipartition $(A,B)$ is a bipartite robust
$(\nu',\tau')$-expander with bipartition $(A,B)$.

We will need a decomposition of dense almost regular graphs due to
K\"uhn, Lo, Osthus, and Staden~\cite{KLOS2015}. To state it, we
recall their definitions. Let $G$ be a graph on $n$ vertices. For
$X\subseteq V(G)$, let $\ol X=V(G)\setminus X$. Let $0<\rho,\nu,\tau<1$.
A set $U\subseteq V(G)$ is a \emph{$\rho$-component} of $G$ if
$|U|\ge\sqrt\rho\,n$ and $e_G(U,\ol U)\le\rho n^2$. For
$U\subseteq V(G)$, the graph $G[U]$ is \emph{$\rho$-close to bipartite
with bipartition $(U_1,U_2)$} if
\begin{enumerate}[label=(C\arabic*)]
\item $U$ is the union of two disjoint sets $U_1$ and $U_2$ with
      $|U_1|,|U_2|\ge\sqrt\rho\,n$;
\item $\bigl||U_1|-|U_2|\bigr|\le\rho n$;
\item $e_G(U_1,\ol{U_2})+e_G(U_2,\ol{U_1})\le\rho n^2$.
\end{enumerate}
The graph $G[U]$ is a \emph{$(\rho,\nu,\tau)$-robust expander
component} of $G$ if
\begin{enumerate}[label=(E\arabic*)]
\item $U$ is a $\rho$-component of $G$;
\item $G[U]$ is a robust $(\nu,\tau)$-expander.
\end{enumerate}
The graph $G[U]$ is a \emph{bipartite $(\rho,\nu,\tau)$-robust
expander component} of $G$ with respect to $(A,B)$ if
\begin{enumerate}[label=(B\arabic*)]
\item $G[U]$ is $\rho$-close to bipartite with bipartition $(A,B)$;
\item $G[U]$ is a bipartite robust $(\nu,\tau)$-expander with
      bipartition $(A,B)$.
\end{enumerate}
Finally, let $D>0$, and let $k$ and $\ell$ be nonnegative integers.
A partition $\cV$ of $V(G)$ is a \emph{robust partition of $G$ with
parameters $\rho,\nu,\tau,k,\ell$}, relative to $D$, if the following
hold.
\begin{enumerate}[label=(D\arabic*)]
\item $\cV=\{V_1,\ldots,V_k,W_1,\ldots,W_\ell\}$ is a partition of
      $V(G)$.
\item For each $i\in[k]$, $G[V_i]$ is a $(\rho,\nu,\tau)$-robust
      expander component of $G$.
\item For each $j\in[\ell]$, there is a partition $(A_j,B_j)$ of
      $W_j$ such that $G[W_j]$ is a bipartite $(\rho,\nu,\tau)$-robust
      expander component of $G$ with respect to $(A_j,B_j)$.
\item For all $X,X'\in\cV$ and all $x\in X$, we have
      $d_G(x,X)\ge d_G(x,X')$.
\item For each $j\in[\ell]$, we have $d_G(u,B_j)\ge d_G(u,A_j)$
      for all $u\in A_j$ and $d_G(v,A_j)\ge d_G(v,B_j)$ for all
      $v\in B_j$.
\item $k+2\ell\le\bigl\lfloor(1+\rho^{1/3})\,n/D\bigr\rfloor$.
\item For all $X\in\cV$, all but at most $\rho n$ vertices $x\in X$
      satisfy $d_G(x,X)\ge D-\rho n$.
\end{enumerate}
In~\cite[Section~3]{KLOS2015}, this definition is made for
$D$-regular graphs $G$. For almost regular graphs it is applied with
$D=\delta(G)$ \cite[Section~5]{KLOS2015}, and the robust partition
theorem then states the following.

\begin{theorem}[K\"uhn, Lo, Osthus, and Staden
{\cite[Theorem~5.1]{KLOS2015}}]\label{thm:robustpartition}
For all $0<\alpha,\tau<1$ and every nondecreasing function
$f:(0,1)\to(0,1)$, there exist an integer $n_0$ and $\gamma>0$ such
that the following holds. For every graph $G$ on $n\ge n_0$ vertices
with
\[
\alpha n\le\delta(G)\le\Delta(G)\le\delta(G)+\gamma n,
\]
there exist reals $\rho,\nu$ with
\[
1/n_0,\ \gamma\le\rho\le\nu\le\tau,\qquad
\rho\le f(\nu),\qquad 1/n_0\le f(\rho),
\]
and nonnegative integers $k,\ell$ such that $G$ has a robust
partition with parameters $\rho,\nu,\tau,k,\ell$, relative to
$D=\delta(G)$.
\end{theorem}

We will also need the following results of K\"uhn, Lo, Osthus, and
Staden, which provide spanning linkages inside the clusters.

\begin{theorem}[K\"uhn, Lo, Osthus, and Staden
\cite{KLOS2015}]\label{thm:linkagetools}
Suppose $0<1/n\ll\nu\ll\tau\ll\eta<1$.
\begin{enumerate}[label=(\roman*)]
\item If $G$ is a robust $(\nu,\tau)$-expander on $n$ vertices with
      $\delta(G)\ge\eta n$, then $G$ is Hamilton $p$-linked for every
      $1\le p\le\nu^4n$.
\item Let $G$ be a bipartite graph with bipartition $(A,B)$, where
      $|A|=|B|=n$, and suppose that $G$ is a bipartite robust
      $(\nu,\tau)$-expander with bipartition $(A,B)$ and
      $\delta(G)\ge\eta n$. If $1\le p\le\nu^4n$, then for any $2p$
      distinct vertices of $G$, exactly $p$ in each part, and any
      pairing of them, $G$ has a spanning linkage for these pairs.
      Individual pairs may have both vertices in the same part.
\end{enumerate}
\end{theorem}

Part (i) is \cite[Corollary~6.9(ii)]{KLOS2015}, and part (ii) is
\cite[Lemma~6.12]{KLOS2015}. Note that the graph in part (ii) has $2n$
vertices, and its robust neighborhoods and expansion are measured with
respect to its order $2n$. We also use the following fact about
deleting a few vertices from a bipartite robust expander.

\begin{lemma}[K\"uhn, Lo, Osthus, and Staden
{\cite[Lemma~4.9(i)]{KLOS2015}}]\label{lem:deletion}
Suppose $0<\nu\ll\tau\ll1$. Let $G$ be a bipartite robust
$(\nu,\tau)$-expander with bipartition $(A,B)$, and let
$W\subseteq V(G)$ with $|W|\le\nu|A|/2$. Then $G-W$ is a bipartite
robust $(\nu/2,2\tau)$-expander with bipartition
$(A\setminus W,B\setminus W)$.
\end{lemma}

Consider $p$ pairs of $2p$ distinct vertices of a bipartite graph with
parts $A$ and $B$, and let $a$ and $b$ be the numbers of pairs inside
$A$ and inside $B$, respectively. If $w$ pairs join the two parts,
then $2a+w$ of the vertices lie in $A$ and $p=a+b+w$. Hence the
condition in Theorem~\ref{thm:linkagetools}(ii) that exactly $p$ of
the vertices lie in each part is equivalent to $a=b$. We prove that,
when the parts have different sizes, the corresponding condition is
$a-b=|A|-|B|$. Bipartite robust expanders whose parts have different
sizes are also treated in~\cite[Section~6]{KLOS2015}, by means of
path systems with as many ends in $A$ as in $B$ for which the numbers
of vertices of $A$ and of $B$ not in the interiors of the paths are
equal. We derive the form that we need from
Theorem~\ref{thm:linkagetools}(ii).

\begin{lemma}\label{lem:unbalanced}
Suppose $0<1/n\ll\nu\ll\tau\ll\eta<1$, and let
$0<\beta\le\nu^4/10^6$. Let $G$ be a bipartite graph on $n$ vertices
with bipartition $(A,B)$ and $\delta(G)\ge\eta n$, and suppose that $G$
is a bipartite robust $(\nu,\tau)$-expander with this bipartition. Let $1\le p\le\beta n$, and let
$x_1,y_1,\ldots,x_p,y_p$ be distinct vertices of $G$. Let $a$ and $b$
be the numbers of indices $i\in[p]$ with $\{x_i,y_i\}\subseteq A$
and with $\{x_i,y_i\}\subseteq B$, respectively. If
\[
a-b=|A|-|B|,
\]
then $G$ has a spanning linkage for the pairs $\{x_i,y_i\}$,
$i\in[p]$.
\end{lemma}

\begin{proof}
We delete $|a-b|$ vertices from the larger part to get a balanced
bipartite graph. As Lemma~\ref{lem:deletion} allows the deleted set to
lie in either part and still yields a bipartite robust expander, we
prove only the case $a\ge b$. By renumbering the pairs, we may assume
that $\{x_i,y_i\}\subseteq B$ for $i\in[b]$, that
$\{x_i,y_i\}\subseteq A$ for $i\in[b+1,a+b]$, and that $x_i$ and $y_i$
lie in distinct parts for $i\in[a+b+1,p]$. Let
$X=\{x_1,y_1,\ldots,x_p,y_p\}$ and $I=[2b+1,a+b]$, so that $|I|=a-b$.
The idea is to delete the vertex $y_i$ for every $i\in I$ and to replace
$y_i$ in its pair by a neighbor $y_i^*$ of $y_i$, which makes the pair
system balanced.

For each $i\in I$ in turn, choose
\[
y_i^*\in N_G(y_i)\setminus\bigl(X\cup\{y_j^*:j\in I,\ j<i\}\bigr).
\]
This is possible, since $y_i$ has at least $\eta n$ neighbors, while
$|X|+(a-b)\le3p\le3\beta n<\eta n$. As $y_i\in A$ and $G$ is bipartite,
$y_i^*\in B$. Thus the vertices $y_i^*$, $i\in I$, are distinct
vertices of $B\setminus X$, and $y_iy_i^*\in E(G)$ for each $i\in I$.

Let $Y=\{y_i:i\in I\}\subseteq A$, let $G^*=G-Y$ with
$A^*=A\setminus Y$, and let $n^*=n-(a-b)$. Then
$|A^*|=|A|-(a-b)=|B|=n^*/2$. By the minimum-degree condition,
$|A|\ge\eta n$, and so, as $\nu\ll\eta$,
\[
|Y|=a-b\le p\le\beta n\le\frac{\nu^4n}{10^6}\le\frac{\nu\eta n}2
\le\frac{\nu|A|}2.
\]
Hence Lemma~\ref{lem:deletion} shows that $G^*$ is a bipartite robust
$(\nu/2,2\tau)$-expander with bipartition $(A^*,B)$. Since every
vertex of $G^*$ loses at most $a-b\le\eta n/2$ neighbors, we have
$\delta(G^*)\ge\eta n/2$. Since $n\ge n^*\ge n-\beta n\ge n/2$, we have
$|A^*|=|B|=n^*/2\ge n/4$, $\delta(G^*)\ge\eta n/2\ge\eta n^*/2$, and
\[
p\le\beta n\le\frac{\nu^4n}{10^6}\le\frac{\nu^4n^*}{32}
=\Bigl(\frac\nu2\Bigr)^4\frac{n^*}2.
\]

We now define $p$ pairs of vertices of $G^*$: the pairs
$\{x_i,y_i\}$ for $i\in[p]\setminus I$, which are unchanged, and
the pairs $\{x_i,y_i^*\}$ for $i\in I$. These $2p$ vertices are
distinct and lie in $G^*$ by the choice of the vertices $y_i^*$ and the
definition of $Y$. Exactly $p$ of them lie in $A^*$ and $p$ in $B$: the
pairs with $i\in[b]$ contribute $2b$ vertices to $B$, the pairs with
$i\in[b+1,2b]$ contribute $2b$ vertices to $A^*$, and every other pair
contributes one vertex to each part. Hence
Theorem~\ref{thm:linkagetools}(ii), applied to $G^*$ with
$n^*/2,\nu/2,2\tau$ in place of $n,\nu,\tau$, respectively, and with
the same $\eta$ (this is possible, since $1/(n^*/2)\ll\nu/2\ll2\tau\ll\eta$
by the hierarchy, as $n^*/2\ge n/4$), gives a spanning
linkage $Q_1^*,\ldots,Q_p^*$ of $G^*$ for these pairs, where $Q_i^*$
joins $x_i$ and $y_i$ for $i\notin I$, and $x_i$ and $y_i^*$ for
$i\in I$.

Finally, let $Q_i=Q_i^*$ for $i\in[p]\setminus I$, and for $i\in I$
let $Q_i$ be obtained from $Q_i^*$ by adding the vertex $y_i$ and the
edge $y_i^*y_i$. Since $y_i\notin V(G^*)$, $Q_i$ is an
$(x_i,y_i)$-path of $G$. The paths $Q_1,\ldots,Q_p$ are
vertex-disjoint, because the paths $Q_i^*$ are vertex-disjoint and the
added vertices $y_i$, $i\in I$, are distinct and lie outside $G^*$,
and they cover $V(G^*)\cup Y=V(G)$. Thus $Q_1,\ldots,Q_p$ form the
desired spanning linkage of $G$.
\end{proof}

\subsection{Linkage between robust expander components}\label{sec:matching}

In the proof of Theorem~\ref{thm:main}, the clusters of the robust
partition will be joined by edges of $G$ that form a matching, and
these edges must connect all the clusters. To describe how a set of edges
connects the clusters, we use the following multigraph. Let $G$ be a
graph, and let $\cU=\{U_1,\ldots,U_q\}$ be a partition of $V(G)$. Let
\[
E_\cU(G)=\bigcup_{i\ne j}E_G(U_i,U_j)
\]
be the set of edges of $G$ joining distinct parts of $\cU$. The
\emph{partition multigraph} $P(G,\cU)$ of $G$ with respect to $\cU$ is
the multigraph with vertex set $\cU$ and edge set $E_\cU(G)$, in which
an edge $uv$ with $u\in U_i$ and $v\in U_j$ joins the vertices $U_i$
and $U_j$. Thus $P(G,\cU)$ is obtained from $G$ by deleting the edges
inside the parts and identifying each part to a single vertex; it has
no loops, but distinct edges of $G$ between the same two parts become
parallel edges. A matching of $G$ contained in $E_\cU(G)$ is called a
\emph{$\cU$-matching} of $G$. For a $\cU$-matching $M$, the
multigraph $P(G,\cU)[M]$ is the subgraph of $P(G,\cU)$ formed by the
edges of $M$; its vertices are the parts that contain an endvertex of
an edge of $M$. An \emph{Euler circuit} of a multigraph is a closed
trail that contains all its edges; we also regard an Euler circuit as
a subgraph. An Euler circuit $T$ of a subgraph of $P(G,\cU)$ is an
\emph{Euler $\cU$-circuit} if $T$ passes through every vertex of
$P(G,\cU)$ and the edges of $G$ corresponding to $E(T)$ form a
$\cU$-matching.

For a multigraph $H$ and a partition $\cP$ of $V(H)$, we say that an
edge of $H$ \emph{crosses} $\cP$ if its two ends lie in distinct sets
of $\cP$. The following theorem of Nash-Williams and Tutte
characterizes the multigraphs that contain $k$ edge-disjoint spanning
trees.

\begin{theorem}[Nash-Williams~\cite{NashWilliams1961},
Tutte~\cite{Tutte1961}]\label{thm:nwt}
Let $k$ be a positive integer. A multigraph $H$ contains $k$
edge-disjoint spanning trees if and only if, for every partition
$\cP$ of $V(H)$ into $s$ sets, at least $k(s-1)$ edges of $H$ cross
$\cP$.
\end{theorem}

We want a $\cU$-matching $M$ of $G$ such that $P(G,\cU)[M]$ contains a
connected spanning even subgraph. It suffices to find $M$ such that
$P(G,\cU)[M]$ contains two edge-disjoint spanning trees of $P(G,\cU)$.
The lemma below provides such an $M$. It can be viewed as a
degree-constrained version of Theorem~\ref{thm:nwt}: under a
toughness-type condition on $G$ and a lower bound on the sizes of the
parts, the
$k$ edge-disjoint spanning trees of $P(G,\cU)$ can be chosen so that
their edges form a matching of $G$. Its proof uses
an iterative exchange argument that goes back to Win~\cite{Win1989},
who proved that, for every integer $d\ge3$, a connected graph $G$ with
$c(G-X)\le(d-2)|X|+2$ for every $X\subseteq V(G)$ has a spanning tree
of maximum degree at most $d$.
Hasanvand~\cite[Theorem~3.1]{Hasanvand2024} extended this approach,
and later refinements of it, from spanning trees to unions of
edge-disjoint spanning trees with degree bounds; we use the same scheme, with maximal tight sets in
place of tree-connected components.

\begin{lemma}\label{lem:twotrees}
	Let $k$ be a positive integer, and let $G$ be a graph such that
	\begin{equation}\label{eq:matching_toughness}
		c(G-S)\le\frac{|S|+2k+1}{2k}
		\qquad\text{for every }S\subseteq V(G).
	\end{equation}
	Let $\cU=\{U_1,\ldots,U_q\}$ be a partition of $V(G)$, where $q\ge2$
	and $|U_i|\ge k(q-1)$ for every $i\in[q]$. Then $G$ has a
	$\cU$-matching $M$ such that $P(G,\cU)[M]$ is the union of $k$
	edge-disjoint spanning trees of $P(G,\cU)$. In particular, if $G$ is
	$4$-tough and $|U_i|\ge2(q-1)$ for every $i\in[q]$, then $P(G,\cU)$
	has an Euler $\cU$-circuit $T$ with at most $2(q-1)$ edges.
\end{lemma}

\begin{proof}
	Our goal is to show the existence of a $\cU$-matching $M$ with exactly
	$k(q-1)$ edges such that, for any partition of $V(P(G,\cU))=\cU$
	into $s$ sets, at least $k(s-1)$ edges of $M$ cross it. For such
	$M$, every part $U_i$ contains an endvertex of an edge of $M$, since
	the partition $\{\{U_i\},\cU\setminus\{U_i\}\}$ is crossed by at
	least $k$ edges of $M$; thus $V(P(G,\cU)[M])=\cU$. By
	Theorem~\ref{thm:nwt} applied to $P(G,\cU)[M]$, this multigraph
	contains $k$ edge-disjoint spanning trees, which are spanning trees
	of $P(G,\cU)$, and since these trees have $k(q-1)=|M|$ edges in
	total, their union is $P(G,\cU)[M]$. To find such an $M$, we use the
	following notation and concepts.
	
	For $I\subseteq[q]$, we let $U_I=\bigcup_{i\in I}U_i$. Let
	$F\subseteq E_\cU(G)$. For $S\subseteq V(G)$, we
	let $F[S]$ be the set of edges of $F$ with both ends contained in
	$S$. We call $F$ \emph{sparse} if
	\begin{equation}\label{eq:matching_sparsity}
		|F[U_I]|\le k(|I|-1)
		\qquad\text{for every nonempty }I\subseteq[q].
	\end{equation}
	A nonempty set $I\subseteq[q]$ is \emph{tight} for $F$ if equality
	holds in \eqref{eq:matching_sparsity}. Every singleton is tight, since
	$F$ contains no edge with both ends in one part. It now suffices to
	show that $G$ has a sparse $\cU$-matching $M$ of size $k(q-1)$. Indeed,
	for such $M$ and any partition $I_1,\ldots,I_s$ of $[q]$, the number of
	edges of $M$ joining distinct sets $U_{I_h}$ is
	\[
	|M|-\sum_{h=1}^s|M[U_{I_h}]|\ge k(q-1)-\sum_{h=1}^sk(|I_h|-1)=k(s-1).
	\]
	
	We first prove some properties of tight sets that will allow us to
	make local modifications of a given matching.
	
	\begin{claim}\label{cl:matching_tightsets}
		Let $F\subseteq E_\cU(G)$ be sparse. The following
		hold.
		\begin{enumerate}[label=(\alph*)]
			\item If $I\subseteq[q]$ and $J\subseteq[q]$ are intersecting tight
			sets for $F$, then $I\cap J$ and $I\cup J$ are also tight for
			$F$. Consequently, the maximal tight sets for $F$ form a
			partition of $[q]$, and every tight set for $F$ is contained in
			a maximal tight set for $F$.
			\item Let $I_1,\ldots,I_s$ be a partition of $[q]$ into tight sets for
			$F$, and let $J\in\{I_1,\ldots,I_s\}$. Let $F'$ be obtained from
			$F$ by replacing the set $F[U_J]$ by a sparse set of $k(|J|-1)$
			edges with both ends in $U_J$, while keeping all other edges of
			$F$. Then $F'$ is sparse, and $I_1,\ldots,I_s$ are tight for $F'$.
			Moreover, if $I_1,\ldots,I_s$ are the maximal tight sets for $F$,
			then they are also the maximal tight sets for $F'$.
			\item Suppose that
			$xy\in E_\cU(G)\setminus F$, and that some
			tight set $I$ for $F$ satisfies $x,y\in U_I$. Then there is a
			unique smallest such tight set $J$, and for every
			$f\in F[U_J]$, the edge set $(F\setminus\{f\})\cup\{xy\}$ is
			sparse.
		\end{enumerate}
	\end{claim}
	
	\begin{claimproof}{cl:matching_tightsets}
		For any $I,J\subseteq[q]$, every edge of $F[U_I]\cup F[U_J]$ has
		both ends in $U_{I\cup J}$, and $F[U_I]\cap F[U_J]=F[U_{I\cap J}]$.
		Hence
		\begin{equation}\label{eq:matching_submodular}
		|F[U_I]|+|F[U_J]|
		=|F[U_I]\cup F[U_J]|+|F[U_{I\cap J}]|
		\le|F[U_{I\cup J}]|+|F[U_{I\cap J}]|.
		\end{equation}
		
		For (a), let $I\cap J\ne\emptyset$, with $I$ and $J$ tight. Since
		$F$ is sparse, we have $|F[U_{I\cup J}]|\le k(|I\cup J|-1)$ and
		$|F[U_{I\cap J}]|\le k(|I\cap J|-1)$. Now since $I$ and $J$ are tight
		for $F$, we have, by \eqref{eq:matching_submodular},
		\[
		\begin{aligned}
			k(|I|-1)+k(|J|-1)
			&=|F[U_I]|+|F[U_J]|\\
			&\le|F[U_{I\cup J}]|+|F[U_{I\cap J}]|\\
			&\le k(|I\cup J|-1)+k(|I\cap J|-1)\\
			&=k(|I|-1)+k(|J|-1).
		\end{aligned}
		\]
		This gives $|F[U_{I\cup J}]|+|F[U_{I\cap J}]|=k(|I|-1)+k(|J|-1)$.
		Combining $|F[U_{I\cup J}]|\le k(|I\cup J|-1)$ and
		$|F[U_{I\cap J}]|\le k(|I\cap J|-1)$, we conclude that $I\cap J$ and
		$I\cup J$ are also tight for $F$. Distinct maximal tight sets are
		therefore disjoint. Since every singleton is tight and lies in a
		maximal tight set, the maximal tight sets form a partition of $[q]$.
		Finally, a tight set $J$ meets some maximal tight set $I$, and
		$J\cup I$ is tight, so $J\subseteq I$ by the maximality of $I$.
		
		For (b), let $I\subseteq[q]$ be nonempty. If $I\cap J=\emptyset$,
		then $F'[U_I]=F[U_I]$. Otherwise, the edges of $F'[U_I]$ with both
		ends in $U_J$ form $F'[U_{I\cap J}]$, which has at most
		$k(|I\cap J|-1)$ edges since the replacement set is sparse, and the
		edges of $F'[U_I]$ with an end outside $U_J$ are exactly the edges of
		$F[U_I]$ with an end outside $U_J$. Since $J$ is tight and $F$ is
		sparse, \eqref{eq:matching_submodular} shows that their number is
		\[
		|F[U_I]|-|F[U_{I\cap J}]|
		\le|F[U_{I\cup J}]|-|F[U_J]|
		\le k(|I\cup J|-1)-k(|J|-1)=k|I\setminus J|.
		\]
		Thus
		$|F'[U_I]|\le k(|I\cap J|-1)+k|I\setminus J|=k(|I|-1)$, proving that
		$F'$ is sparse.
		
		Simply by counting the number of edges of $F'$ within each
		$U_{I_j}$, it follows that each $I_j$ is tight for $F'$. Now suppose
		that $I_1,\ldots,I_s$ are the maximal tight sets for $F$, and let $I'$
		be a maximal tight set for $F'$. By (a) applied to $F'$, the set
		$I'\cup I_j$ is tight for $F'$ whenever $I'$ meets $I_j$, and so
		$I_j\subseteq I'$ by maximality. Hence $I'$ is the union of some of
		the sets $I_1,\ldots,I_s$. Since $F'$ is obtained from $F$ by only
		replacing edges of $F$ with both ends in $U_J$, and the replacement
		preserves the number of edges inside $U_J$, we get
		$|F[U_{I'}]|=|F'[U_{I'}]|=k(|I'|-1)$, showing that $I'$ is also tight
		for $F$. As $I'$ contains some $I_j$, the maximality of $I_j$ for $F$
		gives $I'=I_j$. Thus $I_1,\ldots,I_s$ are also the maximal tight sets
		for $F'$.
		
		For (c), let $J$ be the intersection of all tight sets $I$ for $F$
		with $x,y\in U_I$. Any two such sets intersect, as both contain the
		index of the part containing $x$, so $J$ is tight by a repeated
		application of (a), and $J$ is the unique smallest such set. Let
		$F'=(F\setminus\{f\})\cup\{xy\}$. For any $I\subseteq[q]$, we have
		$|F'[U_I]|\le|F[U_I]|$ unless $x,y\in U_I$, and $|F'[U_I]|\le|F[U_I]|+1$
		always. We only need to consider $I$ with $x,y\in U_I$ for which
		$|F[U_I]|=k(|I|-1)$, i.e., $I$ is tight for $F$. Then $J\subseteq I$,
		and so $f\in F[U_J]\subseteq F[U_I]$. Hence
		$|F'[U_I]|=|F[U_I]|-1+1=k(|I|-1)$. Thus $F'$ is sparse.
	\end{claimproof}
	
	Let $M$ be a sparse $\cU$-matching of $G$ with maximum size. Suppose
	for a contradiction that $|M|<k(q-1)$. Then $[q]$ is not tight for
	$M$. We shall find an integer $s\ge2$ and a set of at most
	$2k(s-1)-2$ vertices, all saturated by $M$, whose deletion leaves at
	least $s$ components. This contradicts \eqref{eq:matching_toughness}.

	For $W\subseteq V(G)$, recall that $\ol W=V(G)\setminus W$; thus
	$M[\ol W]$ is the matching obtained from $M$ by deleting all edges
	covering a vertex of $W$. For $v\in V(G)$, let $I_{\ol W}(v)$ be the
	maximal tight set for $M[\ol W]$ such that $v\in U_{I_{\ol W}(v)}$;
	this is well defined by Claim~\ref{cl:matching_tightsets}(a). If $I$
	is tight for $M[\ol W]$, then, since $M$ is sparse,
	\begin{equation}\label{eq:matching_residualtight}
		k(|I|-1)=|M[U_I\setminus W]|\le|M[U_I]|\le k(|I|-1),
	\end{equation}
	so $I$ is tight for $M$ and $M[U_I]=M[U_I\setminus W]$: no edge of $M$
	covering a vertex of $W$ has both ends in $U_I$. The same argument
	shows that every tight set for $M[\ol{W'}]$ is tight for $M[\ol W]$ whenever
	$W\subseteq W'$, so the partition into maximal tight sets can only be
	refined as $W$ grows. Let $\mathcal A(W,v)$ be the family of all sparse
	$\cU$-matchings $N$ with $|N|=|M|$ such that $N$ is obtained from $M$
	by modifying only edges with both ends in $U_{I_{\ol W}(v)}\setminus W$,
	that is, $M\mathbin{\triangle}N\subseteq E(G[U_{I_{\ol W}(v)}\setminus W])$.
	In particular, $M\in\mathcal A(W,v)$.
	
	Let $S_0=\emptyset$. For each integer $i\ge0$, define
	\begin{equation}\label{eq:matching_closure}
		S_{i+1}=S_i\cup
		\{v\in V(G):v\text{ is saturated by every }N\in\mathcal A(S_i,v)\}.
	\end{equation}
	Since $M\in\mathcal A(S_i,v)$, every vertex of every $S_i$ is saturated
	by $M$. Let $t$ be the smallest index such that $S_t=S_{t+1}$. Since
	$G$ is a finite graph and $S_0\subseteq S_1\subseteq\cdots$, it follows
	that $t$ exists.
	
	Let $I_1,\ldots,I_s$ be the maximal tight sets for $M[\ol{S_t}]$. By
	\eqref{eq:matching_residualtight}, $I_1,\ldots,I_s$ form a partition
	of $[q]$ into tight sets for $M$. Since $[q]$ is not tight for $M$, it
	follows that $s\ge2$. Our goal is to show that for any two distinct
	$I,J\in\{I_1,\ldots,I_s\}$ and any $xy\in E_G(U_I,U_J)\setminus M$, we
	have $S_t\cap\{x,y\}\ne\emptyset$, and then derive a contradiction
	to~\eqref{eq:matching_toughness}.
	
	Suppose otherwise that there exist distinct $I,J\in\{I_1,\ldots,I_s\}$
	and $xy\in E_G(U_I,U_J)\setminus M$ such that
	$S_t\cap\{x,y\}=\emptyset$. Let $i$ be the smallest index such that
	$x$ and $y$ lie in distinct maximal tight sets for $M[\ol{S_i}]$; then
	$0\le i\le t$, and $x,y\notin S_t\supseteq S_{i+1}$. By the choice of
	$i$ and \eqref{eq:matching_closure}, there are $N_x\in\mathcal A(S_i,x)$
	and $N_y\in\mathcal A(S_i,y)$ that leave $x$ and $y$ unsaturated,
	respectively. Let $I_x=I_{\ol{S_i}}(x)$ and $I_y=I_{\ol{S_i}}(y)$, which are
	distinct, and let $M'$ be obtained from $M$ by replacing the edges of
	$M$ within $U_{I_x}$ and $U_{I_y}$ respectively by those of $N_x$ and
	$N_y$. The set $M'$ is a $\cU$-matching: the two replacements occur in
	the disjoint sets $U_{I_x}$ and $U_{I_y}$, and every edge of $M$
	outside $M[U_{I_x}]\cup M[U_{I_y}]$ belongs to both $N_x$ and $N_y$,
	so it shares no vertex with an edge of $N_x[U_{I_x}]$ or of
	$N_y[U_{I_y}]$. Each replacement preserves the number of edges within
	its set, which is $k(|I_x|-1)$ or $k(|I_y|-1)$ as $I_x$ and $I_y$ are
	tight for $M$ by \eqref{eq:matching_residualtight}, and so $|M'|=|M|$.
	Both $x$ and $y$ are unsaturated by $M'$.
	
	Suppose first that $i=0$. Then $I_x$ and $I_y$ are distinct maximal
	tight sets for $M$, and by Claim~\ref{cl:matching_tightsets}(b), $M'$
	is sparse with the same maximal tight sets as $M$. Adding $xy$ to $M'$
	can violate \eqref{eq:matching_sparsity} only for a tight set for $M'$
	containing the indices of the parts of $x$ and $y$; by
	Claim~\ref{cl:matching_tightsets}(a), such a set would lie in a maximal
	tight set for $M'$ containing both indices, which does not exist.
	Then $M'\cup\{xy\}$ is a larger sparse $\cU$-matching, a
	contradiction.
	
	Thus $i\ge1$, and by the choice of $i$, $I_x\cup I_y$ is contained in
	a single maximal tight set $I_0$ for $M[\ol{S_{i-1}}]$. By
	\eqref{eq:matching_residualtight}, $I_0$ is tight for $M$ and no edge
	of $M[U_{I_0}]$ covers a vertex of $S_{i-1}$. By
	Claim~\ref{cl:matching_tightsets}(c), there is a smallest tight set
	$J_0$ for $M$ with $x,y\in U_{J_0}$, and $J_0\subseteq I_0$. If no edge
	of $M[U_{J_0}]$ covered a vertex of $S_i$, then $J_0$ would be tight
	for $M[\ol{S_i}]$, and $x$ and $y$ would lie in one maximal tight set for
	$M[\ol{S_i}]$ by Claim~\ref{cl:matching_tightsets}(a), a contradiction. Thus
	there exists $f=uv\in M[U_{J_0}]$ such that one endvertex of $f$, say
	$u$, is contained in $S_i$. As no edge of $M[U_{I_0}]$ covers a vertex
	of $S_{i-1}$, we have $u\in S_i\setminus S_{i-1}$. By
	\eqref{eq:matching_residualtight}, $f$ joins two distinct maximal
	tight sets for $M[\ol{S_i}]$, and so $f\in M'$, since the replacements
	change only edges inside $U_{I_x}$ and inside $U_{I_y}$. We let
	$M^*=(M'\setminus\{f\})\cup\{xy\}$.
	
	Then $M^*$ is a $\cU$-matching with $|M^*|=|M|$, since $x$ and $y$ are
	unsaturated by $M'$. We claim that $M^*$ is sparse. 
	Note first that $(M\setminus\{f\})\cup\{xy\}$ is
	sparse by Claim~\ref{cl:matching_tightsets}(c). 
	Then as both $f$ and $xy$ join
	distinct maximal tight sets for $M[\ol{S_i}]$, the
	edges inside each such set are unchanged in $(M\setminus\{f\})\cup\{xy\}$,
	and so these sets remain tight for
	$(M\setminus\{f\})\cup\{xy\}$.
	Now, as $M^*$ is obtained from
	$(M\setminus\{f\})\cup\{xy\}$ by replacing the edges inside $U_{I_x}$
	and inside $U_{I_y}$, Claim~\ref{cl:matching_tightsets}(b), applied
	twice, shows that $M^*$ is sparse. Moreover, $M^*$ leaves $u$
	unsaturated: $f$ is the only edge of $M$ covering $u$, the replacements
	avoid $u$, and $u\ne x,y$, as the replacements lie outside $S_i$,
	$x,y\notin S_i$, and $u\in S_i$. Finally, every
	edge of $M\mathbin{\triangle}M^*$ has both ends in
	$U_{I_0}\setminus S_{i-1}$: the replacements lie inside
	$U_{I_x}\setminus S_i$ and $U_{I_y}\setminus S_i$, the edge
	$f\in M[U_{I_0}]$ avoids $S_{i-1}$, and $x,y\in U_{I_0}\setminus S_{i-1}$.
	As $u\in U_{I_0}$, we have $I_{\ol{S_{i-1}}}(u)=I_0$, and so
	$M^*\in\mathcal A(S_{i-1},u)$ leaves $u$ unsaturated. By
	\eqref{eq:matching_closure}, $u\notin S_i$, a contradiction to the
	choice of $u$.
	
	Now since $M$ is sparse and $I_1,\ldots,I_s$ are tight for $M$, the
	number of edges of $M$ joining distinct sets among
	$U_{I_1},\ldots,U_{I_s}$ is
	\[
	|M|-\sum_{h=1}^sk(|I_h|-1)\le k(q-1)-1-k(q-s)=k(s-1)-1.
	\]
	We let $S^*$ be the set of endvertices of these edges of $M$, so that
	$|S^*|\le2k(s-1)-2$. We have $S_t\subseteq S^*$: every vertex of $S_t$
	is saturated by $M$, and since $I_1,\ldots,I_s$ are tight for
	$M[\ol{S_t}]$, by \eqref{eq:matching_residualtight} the edges of $M$
	covering a vertex of $S_t$ join distinct sets $U_{I_h}$. As every edge
	$xy\in E_G(U_I,U_J)\setminus M$ with distinct
	$I,J\in\{I_1,\ldots,I_s\}$ satisfies $\{x,y\}\cap S_t\ne\emptyset$,
	every edge of $G$ joining distinct sets $U_{I_h}$ has an endvertex in
	$S^*$. Moreover, every $U_{I_h}\setminus S^*$ is nonempty, as
	$|U_j|\ge k(q-1)>k(s-1)-1$ and each edge of $M$ covers at most one
	vertex of $U_j$. It follows that $c(G-S^*)\ge s$. However,
	\eqref{eq:matching_toughness} gives
	\[
	c(G-S^*)\le\frac{|S^*|+2k+1}{2k}
	\le\frac{2k(s-1)-2+2k+1}{2k}=s-\frac1{2k},
	\]
	a contradiction. Hence $|M|=k(q-1)$.
	
	Finally, let $k=2$, so that $P(G,\cU)[M]=T_1\cup T_2$ for
	edge-disjoint spanning trees $T_1$ and $T_2$ of $P(G,\cU)$. The
	number of vertices of odd degree in $T_1$ is even, say $2p$, and we
	pair up these vertices as $(x_1,y_1),\ldots,(x_p,y_p)$. For each
	$i\in[p]$, let $P_i$ be the unique $(x_i,y_i)$-path in $T_2$, and let
	$L$ be the spanning subgraph of $T_2$ whose edge set is the symmetric
	difference of $E(P_1),\ldots,E(P_p)$. Taking symmetric differences
	preserves degrees modulo~$2$, and $d_{P_i}(v)$ is odd if and only if
	$v\in\{x_i,y_i\}$. Hence, for every vertex $v$ of $P(G,\cU)$,
	\[
	d_L(v)\equiv\sum_{i=1}^p d_{P_i}(v)\equiv d_{T_1}(v)\pmod 2.
	\]
	As $T_1$ and $L$ are edge-disjoint, every vertex has even degree in
	$T_1\cup L$. Since $T_1\cup L$ contains $T_1$, it is a connected
	spanning even subgraph of $P(G,\cU)[M]$. Hence it has an Euler
	circuit $T$, which passes
	through every vertex of $P(G,\cU)$, and the edges of $G$
	corresponding to $E(T)$ form a $\cU$-matching, as they lie in $M$.
	Thus $T$ is an Euler $\cU$-circuit with at most $|M|=2(q-1)$ edges.
	Since every $4$-tough graph satisfies \eqref{eq:matching_toughness}
	with $k=2$ (if $c(G-S)\ge2$, then $4c(G-S)\le|S|$, and otherwise
	$c(G-S)\le1<(|S|+5)/4$), this proves the last statement.
\end{proof}

\section{Proof of Theorem~\ref{thm:main}}\label{sec:proof}

Let us explain the idea of the proof. Our overall strategy follows
the approach of K\"uhn, Lo, Osthus, and
Staden~\cite[Sections~1.5 and~6]{KLOS2015}, who reduce the Hamiltonicity of
a graph with a robust partition to finding a system of paths that
connects the clusters and is balanced in the bipartite clusters. The two
tasks that they carry out using $3$-connectivity and regularity,
connecting the clusters and balancing the bipartite clusters, are here
carried out using toughness.

Two difficulties arise in our setting that do not occur in theirs,
and the tools for overcoming them are the main new ingredients of
this paper. First, since the minimum degree may be any linear
function of $n$, the number of clusters is bounded only in terms of
$\alpha$, so the clusters cannot be connected by a case analysis.
This is where toughness first enters, through
Lemma~\ref{lem:twotrees}.
Second, without regularity, the edges needed to
balance the two parts of a bipartite robust expander component need
not lie in that component. We use an auxiliary digraph and the
toughness condition to overcome this challenge.

\begin{proof}
		\setcounter{claim}{0}\stepcounter{claimblock}%
		Since a simple graph on $n$ vertices has maximum degree at most
		$n-1$, the theorem is vacuous when $\alpha\ge1$. Thus let
		$0<\alpha<1$, and let
		\[
		q_0=\lceil2/\alpha\rceil\qquad\text{and}\qquad\eta=\alpha/(4q_0).
		\]
		Below, $q_0$ will be an upper bound on the number of clusters of the
		robust partition of $V(G)$, where each bipartite cluster is counted
		twice, and $\eta n$ will be a lower bound on the minimum degrees inside
		the clusters. We also choose two constants $\tau$ and $\tau_0$ ourselves.
		The constant $\tau$ is the expansion parameter in
		Theorem~\ref{thm:linkagetools}(i) and Lemma~\ref{lem:unbalanced}. We
		will apply Theorem~\ref{thm:robustpartition} with $\tau_0$ in place
		of $\tau$, and so the parameter $\nu$ it provides satisfies
		$\nu\le\tau_0$. We choose $\tau>0$ sufficiently small with respect to
		$\eta$, and then $\tau_0\in(0,\tau]$ sufficiently small with respect
		to $\tau$, so that for every $\mu\in(0,\tau_0]$, in particular for
		$\mu=\nu$, the hierarchy $\mu\ll\tau\ll\eta<1$ required in
		Theorem~\ref{thm:linkagetools}(i) and Lemma~\ref{lem:unbalanced},
		with $\mu$ in place of $\nu$, holds. It remains to make the condition
		$1/n\ll\mu$ in these results explicit. By the convention on
		hierarchies, there is a nondecreasing function
		$g_0\colon(0,1)\to(0,1)$, depending only on $\tau_0$, $\tau$, and
		$\eta$, such that for every $\mu\in(0,\tau_0]$,
		Theorem~\ref{thm:linkagetools}(i) and Lemma~\ref{lem:unbalanced},
		with $\mu$ in place of $\nu$, hold for every graph on at least
		$1/g_0(\mu)$ vertices.

		The parameters $\rho$ and $\nu$ of the robust partition are provided
		by Theorem~\ref{thm:robustpartition} together with the partition, and
		they depend on $G$, so we cannot choose them. The proof needs $\rho$
		to be small, and $n$ to be large, compared with certain functions of
		$\nu$. The only way to guarantee this is to prescribe these functions
		through $f$, where $f\colon(0,1)\to(0,1)$ may be any nondecreasing
		function, as allowed in Theorem~\ref{thm:robustpartition}. Thus $f$
		has to be fixed before Theorem~\ref{thm:robustpartition} is applied;
		the constants $n_0$ and $\gamma$ provided by the theorem depend on
		$f$. Accordingly, we define $f$ by
		\[
		f(\mu)=\min\Bigl\{\eta\, g_0(\mu),\ \frac{\eta}{4q_0+1},\
		\frac{\eta}{100q_0},\ \frac{\mu^4\eta}{16\cdot10^8q_0^3(q_0+1)}\Bigr\}.
		\]
		Then $f$ is nondecreasing and depends only on $\alpha$, $\tau$, and
		$\tau_0$, and hence only on $\alpha$.
		
		We apply Theorem~\ref{thm:robustpartition} with $\alpha$, with our
		constant $\tau_0$ in place of $\tau$, and with $f$, and let $\varepsilon=\gamma$
		and $n_0$ be the resulting constants; they depend only on $\alpha$.
		Let $G$ be a $4$-tough graph on $n\ge n_0$ vertices satisfying
		\eqref{eq:main}. We show that $G$ is Hamiltonian.
		
		\medskip
		\noindent\textbf{The robust partition.}
		By Theorem~\ref{thm:robustpartition}, $G$ has a robust partition
		$\cV=\{V_1,\ldots,V_k,W_1,\ldots,W_\ell\}$ with parameters
		$\rho,\nu,\tau_0,k,\ell$, relative to $D=\delta(G)$, where
		\begin{equation}\label{eq:partitionparameters}
			1/n_0,\ \varepsilon\le\rho\le\nu\le\tau_0\qquad\text{and}\qquad
			\rho\le f(\nu).
		\end{equation}
		Let
		\[
		\beta=\frac{\nu^4}{10^6}\qquad\text{and}\qquad
		\lambda=\frac{\beta\eta}{100q_0(q_0+1)},
		\]
		so that the last term in the definition of $f(\nu)$ is
		$\lambda/(16q_0^2)$. By \eqref{eq:partitionparameters}, we have
		$1/n\le1/n_0\le\rho\le f(\nu)$. Hence the definition of $f$ gives
		\begin{equation}\label{eq:conditions}
			\eta n\ge\frac1{g_0(\nu)},\qquad \eta n\ge4q_0+1,\qquad
			q_0\rho\le\frac{\eta}{100},\qquad 16q_0^2\le\lambda n,\qquad
			8q_0^2\rho\le\lambda,
		\end{equation}
		where the first, second, and fourth inequalities follow from
		$1/n\le f(\nu)$, and the third and fifth from $\rho\le f(\nu)$. We
		need the first inequality to apply Theorem~\ref{thm:linkagetools}(i)
		and Lemma~\ref{lem:unbalanced} inside the clusters, the second to
		ensure $n\ge3$ and to apply Lemma~\ref{lem:twotrees}, and the third
		and the last two later on, respectively, in the proof of
		Claim~\ref{cl:sink} and in \eqref{eq:budget}.
		
		Let $q=k+\ell$, and for each $j\in[\ell]$, let $(A_j,B_j)$ be the
		bipartition of $W_j$ provided by property (D3). Property (D6),
		property (D4), and properties (D4) and (D5) together give,
		respectively,
		\[
		k+2\ell\le(1+\rho^{1/3})\frac{n}{\delta(G)},\qquad
		\delta(G[V_i])\ge\frac{\delta(G)}q,\qquad
		\delta(G[A_j,B_j])\ge\frac{\delta(G)}{2q}.
		\]
		Since $\rho<1$ and $\delta(G)\ge\alpha n$, we have
		$q\le k+2\ell<2/\alpha\le q_0$. Hence the last two bounds are at
		least $\alpha n/q_0=4\eta n$ and $\alpha n/(2q_0)=2\eta n$,
		respectively. Property (D2) states that each $G[V_i]$ is a
		$(\rho,\nu,\tau_0)$-robust expander component of $G$, and so, by
		(E2), $G[V_i]$ is a robust $(\nu,\tau_0)$-expander. Property (D3)
		states that each $G[W_j]$ is a bipartite $(\rho,\nu,\tau_0)$-robust
		expander component with respect to $(A_j,B_j)$. By (B2), the graph
		$G[W_j]$ is a bipartite robust $(\nu,\tau_0)$-expander with the
		ordered bipartition $(A_j,B_j)$, and by (B1) it is $\rho$-close to
		bipartite with this bipartition, so (C2) gives
		$\bigl||A_j|-|B_j|\bigr|\le\rho n$. Since deleting the edges inside
		$A_j$ and inside $B_j$ does not affect the bipartite expansion
		condition, $G[A_j,B_j]$ is also a bipartite robust
		$(\nu,\tau_0)$-expander with bipartition $(A_j,B_j)$. Finally, as
		$\tau_0\le\tau$, every robust $(\nu,\tau_0)$-expander is a robust
		$(\nu,\tau)$-expander, and likewise in the bipartite case.
		
		Recall that the members of $\cV$ are called clusters. Each $V_i$ is an
		\emph{ordinary} cluster, and each $W_j=A_j\cup B_j$ is a
		\emph{bipartite} cluster with \emph{parts} $A_j$ and $B_j$. We shall use
		only the following properties of the partition:
		\begin{enumerate}[label=(\roman*)]
			\item $1\le q\le k+2\ell\le q_0$;
			\item for each $i\in[k]$, $G[V_i]$ is a robust $(\nu,\tau)$-expander
			with $\delta(G[V_i])\ge\eta n$;
			\item for each $j\in[\ell]$, $G[A_j,B_j]$ is a bipartite robust
			$(\nu,\tau)$-expander with bipartition $(A_j,B_j)$,
			$\delta(G[A_j,B_j])\ge\eta n$, and
			$\bigl||A_j|-|B_j|\bigr|\le\rho n$.
		\end{enumerate}
		By (ii) and (iii), every cluster, and every part of a bipartite cluster,
		has at least $\eta n$ vertices.
		
		\medskip
		\noindent\textbf{Spanning linkages inside the clusters.}
		By \eqref{eq:conditions}, we have $n\ge3$, and every cluster $V$ has
		$|V|\ge\eta n\ge1/g_0(\nu)$. Moreover, $\nu\le\tau_0$, the minimum degree of
		$G[V_i]$ or of $G[A_j,B_j]$ is at least $\eta n\ge\eta|V|$, and these
		graphs are robust, or bipartite robust, $(\nu,\tau)$-expanders by
		(ii) and (iii). Hence Theorem~\ref{thm:linkagetools}(i) and
		Lemma~\ref{lem:unbalanced}, both applied with $s$ in place of $p$,
		the latter also with $\beta=\nu^4/10^6$, give the following.
		\begin{equation}\label{eq:locallinkage}
			\begin{minipage}{0.88\textwidth}
				Let $V\in\cV$, let $1\le s\le\beta|V|$, and let
				$\{w_1,w_1'\},\ldots,\{w_s,w_s'\}$ be pairs of $2s$ distinct vertices
				of $V$. If $V=V_i$, then $G[V_i]$ has a spanning linkage for these
				pairs. If $V=W_j$ and the number of these pairs inside $A_j$ minus the
				number inside $B_j$ equals $|A_j|-|B_j|$, then $G[A_j,B_j]$ has a
				spanning linkage for these pairs.
			\end{minipage}
		\end{equation}
		
		\medskip
		\noindent\textbf{Alternating cycles.}
		Let $G^+$ be the multigraph obtained from $G$ by adding, for any two
		distinct vertices $u,v$ of the same cluster, a new edge $uv$, which we
		call an \emph{auxiliary edge}. An \emph{alternating cycle} is a cycle
		of $G^+$, possibly of length two, whose edges are alternately edges of
		$G$ and auxiliary edges. Our goal is to find an alternating cycle to
		which the following claim applies.
		
		\begin{claim}\label{cl:skeleton}
			Let $\Gamma$ be an alternating cycle such that
			\begin{enumerate}[label=(\alph*)]
				\item every cluster contains at least one auxiliary edge of $\Gamma$;
				\item every cluster $V$ contains at most $\beta|V|$ auxiliary edges of
				$\Gamma$;
				\item for every $j\in[\ell]$, the number of auxiliary edges of
				$\Gamma$ inside $A_j$ minus the number inside $B_j$ equals
				$|A_j|-|B_j|$.
			\end{enumerate}
			Then $G$ is Hamiltonian.
		\end{claim}
		
		\begin{claimproof}{cl:skeleton}
			Let $V$ be a cluster. As the vertices of $\Gamma$ are distinct, the
			auxiliary edges of $\Gamma$ inside $V$ are $s$ pairs of $2s$ distinct
			vertices of $V$, where $1\le s\le\beta|V|$ by (a) and (b). By (c) and
			\eqref{eq:locallinkage}, the graph $G[V]$, or $G[A_j,B_j]$ if
			$V=W_j$, has a spanning linkage for these pairs. Replacing every
			auxiliary edge of $\Gamma$ by the path of these linkages joining its
			ends turns $\Gamma$ into a closed walk in $G$, in which these paths
			alternate with the edges of $G$ on $\Gamma$: each edge $xy$ of $G$ on
			$\Gamma$ joins the path replacing the auxiliary edge of $\Gamma$ at
			$x$ to the path replacing the auxiliary edge of $\Gamma$ at $y$. Since the clusters
			partition $V(G)$ and the vertices of $\Gamma$ are exactly the ends of
			these paths, this closed walk passes through every vertex of $G$
			exactly once. As $n\ge3$, it is a Hamiltonian cycle.
		\end{claimproof}
		
		\medskip
		\noindent\textbf{The Euler circuit $T_0$.}
		We first find an alternating cycle $\Gamma_0$ with $m_0\le2q_0$
		auxiliary edges that satisfies (a). If $q=1$, let $xy$ be an edge of
		$G$, let $\Gamma_0$ be the alternating cycle of length two formed by
		the edge $xy$ and the auxiliary edge $xy$, and let $m_0=1$.
		
		So let $q\ge2$. Every cluster has at least $\eta n\ge4q_0+1>2(q-1)$
		vertices. Since $G$ is $4$-tough, Lemma~\ref{lem:twotrees} gives an
		Euler $\cV$-circuit $T_0$ with at most $2(q-1)$ edges. Let $M_0$ be
		the $\cV$-matching of $G$ corresponding to $E(T_0)$, and let
		$m_0=|M_0|$, so that $m_0\le2(q-1)\le2q_0$. As $T_0$ passes through
		$q\ge2$ clusters, $m_0\ge2$. We regard $T_0$ cyclically. At each visit of $T_0$
		to a cluster $V$, we pair the endvertex in $V$ of the edge of $M_0$ by
		which $T_0$ enters $V$ with the endvertex in $V$ of the edge by which
		it leaves $V$; we call the resulting pairs the \emph{pairs of $T_0$}.
		Since $T_0$ traverses every edge of $M_0$ once and $M_0$ is a matching,
		numbering the visits in the order of $T_0$ gives the pairs of $T_0$
		as $\{x_i,y_i\}$, $i\in[m_0]$, where $x_1,y_1,\ldots,x_{m_0},y_{m_0}$ are
		distinct, $x_i$ and $y_i$ lie in the cluster of the $i$-th visit, and
		$M_0=\{y_ix_{i+1}:i\in[m_0]\}$ with $x_{m_0+1}=x_1$. Hence
		$\Gamma_0=x_1y_1x_2y_2\cdots x_{m_0}y_{m_0}x_1$, whose auxiliary edges are
		the pairs of $T_0$ and whose edges of $G$ are those of $M_0$, is an
		alternating cycle. As $T_0$ visits every cluster, $\Gamma_0$ satisfies
		(a).
		
		For $j\in[\ell]$, let $p_{A_j}$ and $p_{B_j}$ be the numbers of
		auxiliary edges of $\Gamma_0$ inside $A_j$ and inside $B_j$,
		respectively, and let
		\[
		d_j=\bigl(|A_j|-p_{A_j}\bigr)-\bigl(|B_j|-p_{B_j}\bigr).
		\]
		Then $\Gamma_0$ satisfies (c) if and only if $d_j=0$ for every $j$.
		Interchanging the names of $A_j$ and $B_j$ replaces $d_j$ by $-d_j$.
		Since (iii) is used below only through \eqref{eq:locallinkage}, which
		is symmetric under this interchange, and through the symmetric bounds
		$|A_j|,|B_j|\ge\eta n$ and $\bigl||A_j|-|B_j|\bigr|\le\rho n$, we may
		assume that $d_j\ge0$ for every $j$. We now modify $\Gamma_0$ into an
		alternating cycle $\Gamma$ so that (c) is satisfied. To this end, for
		every $j$, we construct $d_j$ alternating cycles, each with exactly one
		auxiliary edge inside $A_j$ and no other auxiliary edge inside a part
		of a bipartite cluster, and merge them into $\Gamma_0$. These cycles
		are found by means of the following digraph.
		
		\medskip
		\noindent\textbf{The type digraph.}
		Each ordinary cluster $V_i$ is a \emph{type}, and the two parts $A_j,B_j$
		of each bipartite cluster $W_j$ are two types. Thus the set $\cT$ of all
		types is a partition of $V(G)$ into sets of at least $\eta n$
		vertices, and $|\cT|=k+2\ell\le q_0$. For $z\in V(G)$, let $T(z)$ be
		the type containing $z$. We define an involution $\sigma$ on $\cT$ by
		\[
		\sigma(V_i)=V_i,\qquad \sigma(A_j)=B_j,\qquad \sigma(B_j)=A_j.
		\]
		Recall that $\lambda=\beta\eta/(100q_0(q_0+1))$. An unordered pair
		$\{X,Y\}$ of types, where $X=Y$ is allowed, is \emph{rich} if $G$ has
		a matching of size $\lceil\lambda n\rceil$ with edges from $E_G(X,Y)$; for
		$X=Y$, these are the edges of $G[X]$. For every rich pair $\{X,Y\}$,
		we fix such a matching $M_{X,Y}$. Let $\cD$ be the digraph on $\cT$,
		where loops are allowed, in which $X\to Y$ is an arc if and only if
		$\{\sigma(X),Y\}$ is rich. As $\{\sigma(X),Y\}=\{\sigma(\sigma(Y)),\sigma(X)\}$,
		\begin{equation}\label{eq:reverse}
			X\longrightarrow Y\quad\Longleftrightarrow\quad
			\sigma(Y)\longrightarrow\sigma(X),
		\end{equation}
		and so reversing a directed path of $\cD$ and applying the involution
		to each of its vertices gives again a directed path.
		
		For every pair $\{X,Y\}$ of types that is not rich, we choose a
		maximum matching of $G$ with edges from $E_G(X,Y)$, and we
		let $Z$ be the set of endvertices of all these matchings.
		
		\begin{claim}\label{cl:Z}
			We have $|Z|\le\eta n/100$, and every edge $ww'$ of $G$ such that
			$\{T(w),T(w')\}$ is not rich has an endvertex in $Z$.
		\end{claim}
		
		\begin{claimproof}{cl:Z}
			There are $(k+2\ell)(k+2\ell+1)/2$ unordered pairs of types,
			including the pairs $\{X,X\}$, and each matching
			chosen for a pair $\{X,Y\}$ that is not rich has at most
			$\lceil\lambda n\rceil-1<\lambda n$ edges. Hence, as $k+2\ell\le q_0$ and
			$\beta\le1$,
			\[
			|Z|\le(k+2\ell)(k+2\ell+1)(\lceil\lambda n\rceil-1)<q_0(q_0+1)\lambda n\le\frac{\eta n}{100}.
			\]
			If an edge $ww'$ of $G$ with $\{T(w),T(w')\}$ not rich had no
			endvertex in $Z$, it could be added to the maximum matching chosen
			for $\{T(w),T(w')\}$.
		\end{claimproof}
		
		\begin{claim}\label{cl:sink}
			Let $\cH$ be a strong component of $\cD$ such that no arc of $\cD$
			leaves $\cH$, and let $\sigma(\cH)=\{\sigma(X):X\in V(\cH)\}$. Then
			$\sigma(\cH)=V(\cH)$.
		\end{claim}
		
		\begin{claimproof}{cl:sink}
			By \eqref{eq:reverse}, the involution maps strong components of $\cD$
			to strong components, so $\sigma(\cH)$ is the vertex set of a strong
			component of $\cD$. Suppose to the contrary that $\sigma(\cH)\ne
			V(\cH)$. Then $V(\cH)\cap\sigma(\cH)=\emptyset$, so $V(\cH)$
			contains no ordinary type and at most one part of each bipartite
			cluster. Let
			\[
			W=\bigcup_{X\in V(\cH)}X
			\qquad\text{and}\qquad
			W'=\bigcup_{X\in V(\cH)}\sigma(X).
			\]
			Then $W\cap W'=\emptyset$, $|W'|\ge\eta n$, and, summing
			$\bigl||A_j|-|B_j|\bigr|\le\rho n$ over at most $q_0$ bipartite clusters,
			$|W|\le|W'|+q_0\rho n$. Let $x\in W'\setminus Z$, and suppose that $x$
			has a neighbor $w$ in $G-(W\cup Z)$. As $xw$ has no endvertex in $Z$,
			the pair $\{T(x),T(w)\}$ is rich by Claim~\ref{cl:Z}, that is,
			$\sigma(T(x))\to T(w)$ is an arc of $\cD$. Since $\sigma(T(x))\in V(\cH)$
			and no arc leaves $\cH$, we get $T(w)\in V(\cH)$, that is, $w\in W$, a
			contradiction. Thus every vertex of $W'\setminus Z$ is an isolated
			vertex of $G-(W\cup Z)$, and $|W'\setminus Z|\ge\eta n-\eta n/100\ge2$
			by Claim~\ref{cl:Z}. Since $G$ is $4$-tough,
			\[
			4(|W'|-|Z|)\le|W\cup Z|\le|W'|+q_0\rho n+|Z|.
			\]
			By \eqref{eq:conditions}, we have $q_0\rho n\le\eta n/100$, and so,
			by Claim~\ref{cl:Z}, $3\eta n\le3|W'|\le q_0\rho n+5|Z|\le6\eta n/100$,
			a contradiction.
		\end{claimproof}
		
		\begin{claim}\label{cl:weak}
			Every weak component of $\cD$ is strongly connected and closed under
			the involution.
		\end{claim}
		
		\begin{claimproof}{cl:weak}
			Let $\cD^*$ be a weak component of $\cD$. Since the strong components
			of $\cD^*$ form an acyclic digraph, there is a strong component $\cH$
			of $\cD^*$ such that no arc of $\cD^*$ leaves $\cH$. As $\cD^*$ is a
			weak component, no arc of $\cD$ leaves $\cH$, and so
			$\sigma(\cH)=V(\cH)$ by Claim~\ref{cl:sink}. Suppose that an arc
			$Y\to X$ of $\cD$ enters $\cH$, that is, $X\in V(\cH)$ and
			$Y\notin V(\cH)$. By \eqref{eq:reverse}, $\sigma(X)\to\sigma(Y)$ is an
			arc of $\cD$, where $\sigma(X)\in\sigma(\cH)=V(\cH)$ and
			$\sigma(Y)\notin\sigma(\cH)=V(\cH)$. This arc leaves $\cH$, a
			contradiction. Thus no arc of $\cD$ joins $V(\cH)$ and
			$\cT\setminus V(\cH)$. Since $\cD^*$ is weakly connected, $\cD^*=\cH$.
			Hence $\cD^*$ is strongly connected, and it is closed under the
			involution as $\sigma(\cH)=V(\cH)$.
		\end{claimproof}
		
		For every $j\in[\ell]$, let $\cD^*$ be the weak
		component of $\cD$ containing $B_j$. By Claim~\ref{cl:weak}, $\cD^*$
		is strongly connected and contains $A_j=\sigma(B_j)$, so $\cD$
		contains a directed path $\bP_j$ from $B_j$ to $A_j$, and
		$1\le e(\bP_j)\le k+2\ell-1$.
		
		\medskip
		\noindent\textbf{Balancing $\Gamma_0$.}
		Let $h=\sum_jd_je(\bP_j)$. Since $e(\bP_j)\le k+2\ell-1$,
		$d_j\le\bigl||A_j|-|B_j|\bigr|+p_{B_j}\le\rho n+p_{B_j}$, $\ell\le q_0$,
		and $\sum_jp_{B_j}\le m_0$, we have
		\[
		h\le(k+2\ell-1)\sum_jd_j\le(k+2\ell-1)\sum_j(\rho n+p_{B_j})\le(k+2\ell-1)(q_0\rho n+m_0).
		\]
		By \eqref{eq:conditions}, we have
		$16q_0^2\le\lambda n$ and $8q_0^2\rho\le\lambda$. Since $k+2\ell\le q_0$ and
		$m_0\le2q_0$, we obtain
		\begin{equation}\label{eq:budget}
			m_0+h\le(k+2\ell)m_0+(k+2\ell-1)q_0\rho n\le2q_0^2+q_0^2\rho n\le\frac{\lambda n}4.
		\end{equation}
		
		For every $j\in[\ell]$, we now choose $d_j$ alternating cycles
		$R_{j,1},\ldots,R_{j,d_j}$, one after another, such that all these
		cycles and $\Gamma_0$ are pairwise vertex-disjoint. Suppose that we
		are to choose $R_{j,i}$, write $\bP_j=X_0X_1\cdots X_m$ with $m=e(\bP_j)$,
		$X_0=B_j$, and $X_m=A_j$, and let $L$ be the set of vertices of
		$\Gamma_0$ and of the cycles chosen before. For $r=1,\ldots,m$ in
		turn, the pair $\{\sigma(X_{r-1}),X_r\}$ is rich, as $X_{r-1}\to X_r$ is
		an arc of $\cD$, and we choose an edge $a_rb_r$ of
		$M_{\sigma(X_{r-1}),X_r}$ with $a_r\in\sigma(X_{r-1})$ and $b_r\in X_r$ that
		avoids $L$ and $a_1,b_1,\ldots,a_{r-1},b_{r-1}$; if
		$\sigma(X_{r-1})=X_r$, either end may serve as $a_r$. This is possible,
		since by \eqref{eq:budget} fewer than $2(m_0+h)\le\lambda n/2<\lceil\lambda n\rceil$
		vertices are to be avoided, and each of them lies in at most one edge
		of the matching. For $r\in[m-1]$, the vertices $b_r\in X_r$ and
		$a_{r+1}\in\sigma(X_r)$ lie in one cluster, and in distinct parts of it if
		this cluster is bipartite, while $b_m\in X_m=A_j$ and
		$a_1\in\sigma(X_0)=A_j$. Hence
		\[
		R_{j,i}=a_1b_1a_2b_2\cdots a_mb_ma_1,
		\]
		with the edges $a_rb_r$ of $G$ and the auxiliary edges $b_ra_{r+1}$,
		$r\in[m-1]$, and $b_ma_1$, is an alternating cycle with exactly one
		auxiliary edge, $b_ma_1$, inside a part of a bipartite cluster, namely
		inside $A_j$; see Figure~\ref{fig:route-pairs}. Consequently,
		\begin{equation}\label{eq:Rcount}
			\begin{minipage}{0.88\textwidth}
				for every $j'\in[\ell]$, the number of auxiliary edges of $R_{j,i}$
				inside $A_{j'}$ minus the number inside $B_{j'}$ is $1$ if $j'=j$ and
				$0$ otherwise.
			\end{minipage}
		\end{equation}
		
		\begin{figure}[htbp]
			\centering
			\begin{tikzpicture}[rt picture]
				\node at (-.95,1.5) {$\mathbf P_j:$};
				\node (rtXzero) at (.3,1.5) {$B_j$};
				\node (rtXone) at (3,1.5) {$X_1$};
				\node (rtXtwo) at (6.6,1.5) {$X_2$};
				\node (rtXthree) at (9.3,1.5) {$A_j$};
				\draw[rt map] (rtXzero) -- (rtXone);
				\draw[rt map] (rtXone) -- (rtXtwo);
				\draw[rt map] (rtXtwo) -- (rtXthree);
				
				\path[rt shade] (2.05,-.78) rectangle (3.95,.5);
				\path[rt shade] (5.65,-.78) rectangle (7.55,.5);
				
				\node[rt vertex,label=above:{$a_1$}] (rtaone) at (0,0) {};
				\node[rt vertex,label=above:{$b_1$}] (rtbone) at (2.4,0) {};
				\node[rt vertex,label=above:{$a_2$}] (rtatwo) at (3.6,0) {};
				\node[rt vertex,label=above:{$b_2$}] (rtbtwo) at (6,0) {};
				\node[rt vertex,label=above:{$a_3$}] (rtathree) at (7.2,0) {};
				\node[rt vertex,label=above:{$b_3$}] (rtbthree) at (9.6,0) {};
				\draw[rt represent] (rtaone) -- (rtbone);
				\draw[rt formal] (rtbone) -- (rtatwo);
				\draw[rt represent] (rtatwo) -- (rtbtwo);
				\draw[rt formal] (rtbtwo) -- (rtathree);
				\draw[rt represent] (rtathree) -- (rtbthree);
				\draw[rt formal] (rtbthree.south) .. controls (9.6,-1.45) and (0,-1.45) .. (rtaone.south);
				
				\node at (-.75,0) {$A_j$};
				\node at (2.4,-.47) {$X_1$};
				\node at (3.6,-.47) {$\sigma(X_1)$};
				\node at (6,-.47) {$X_2$};
				\node at (7.2,-.47) {$\sigma(X_2)$};
				\node at (10.35,0) {$A_j$};
				\node[rt note] at (4.8,-1.32) {auxiliary edge inside $A_j$};
				
				\draw[rt formal] (2.1,-1.95) -- (2.85,-1.95);
				\node[rt note,anchor=west] at (2.97,-1.95) {auxiliary edge};
				\draw[rt represent] (5.6,-1.95) -- (6.35,-1.95);
				\node[rt note,anchor=west] at (6.47,-1.95) {edge of $G$};
			\end{tikzpicture}
			\caption{The alternating cycle $R_{j,i}$ for a directed path
				$\mathbf P_j=B_jX_1X_2A_j$ with three arcs. The labels next to the
				vertices give the types containing them. The shaded boxes group the
				auxiliary edges $b_1a_2$ and $b_2a_3$, whose ends lie in one cluster,
				and in distinct parts of it if this cluster is bipartite; the auxiliary
				edge $b_3a_1$ lies inside $A_j$. The edges $a_1b_1$, $a_2b_2$, and
				$a_3b_3$ belong to the matchings $M_{A_j,X_1}$,
				$M_{\sigma(X_1),X_2}$, and $M_{\sigma(X_2),A_j}$,
				respectively.}
			\label{fig:route-pairs}
		\end{figure}
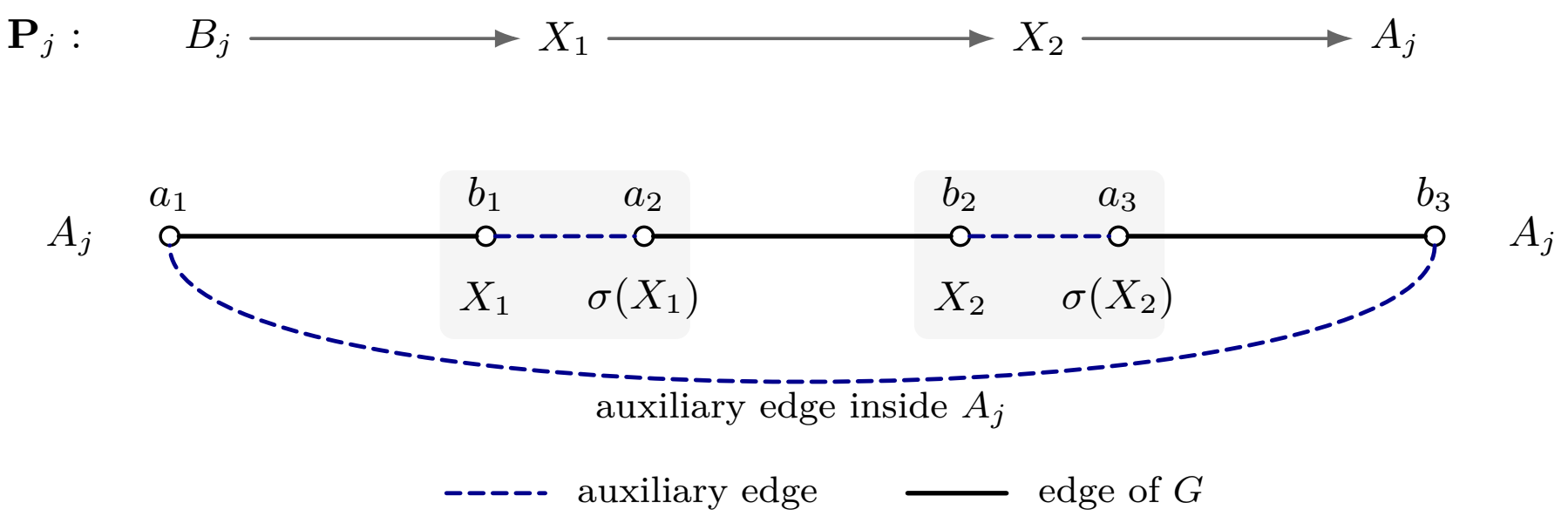
		
		\medskip
		\noindent\textbf{Merging.}
		Let $\Gamma=\Gamma_0$. We merge the cycles $R_{j,i}$ into $\Gamma$ one
		by one, keeping $\Gamma$ an alternating cycle that satisfies (a). Let
		$R=R_{j,i}$ be the next cycle, let $zz'$ be its auxiliary edge inside
		$A_j$, and let $uv$ be an auxiliary edge of $\Gamma$ inside $W_j$,
		which exists by (a). As $V(\Gamma)$ is the union of $V(\Gamma_0)$ and
		the vertex sets of the cycles merged before, $R$ and $\Gamma$ are
		vertex-disjoint, so $u,v,z,z'$ are distinct. Deleting the auxiliary
		edge $uv$ from $\Gamma$ leaves a path of $G^+$ from $u$ to $v$, and
		deleting the auxiliary edge $zz'$ from $R$ leaves a path of $G^+$
		from $z$ to $z'$; this also holds if $\Gamma$ or $R$ has length two.
		Both paths start and end with edges of $G$, so adding the auxiliary
		edges $uz$ and $vz'$ joins them into an alternating cycle with vertex
		set $V(\Gamma)\cup V(R)$; see Figure~\ref{fig:splice}. It contains
		$uz$, which lies inside $W_j$, instead of $uv$, and so it satisfies
		(a). We let $\Gamma$ be this cycle.
		
		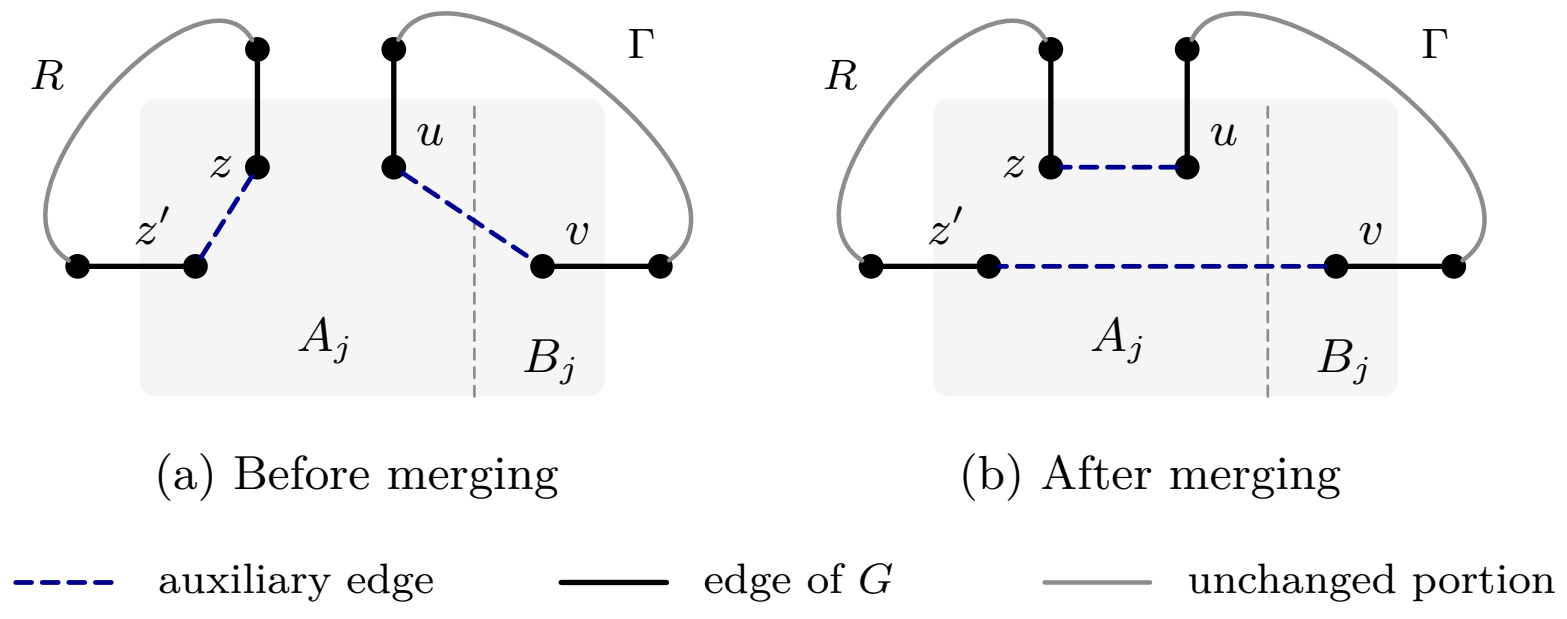
\begin{figure}[htbp]
			\centering
			\begin{tikzpicture}[rt picture,x=.82cm,y=.82cm]
				\foreach \shift/\sfx in {0/A,6.4/B}{
					\begin{scope}[xshift=\shift*.82 cm]
						\path[rt shade] (-.2,.7) rectangle (3.55,3.1);
						\draw[black!45,line width=.5pt,dash pattern=on 2pt off 1.5pt]
						(2.5,.7) -- (2.5,3.1);
						\node at (1.3,1.15) {$A_j$};
						\node at (3.12,.98) {$B_j$};
						\node[rt actual,label={[label distance=-1pt]above right:{$u$}}]
						(rta\sfx) at (1.85,2.55) {};
						\node[rt actual,label={[label distance=-1pt]above right:{$v$}}]
						(rtb\sfx) at (3.05,1.75) {};
						\node[rt actual,label={[label distance=-1pt]left:{$z$}}]
						(rtc\sfx) at (.75,2.55) {};
						\node[rt actual,label={[label distance=-1pt]above left:{$z'$}}]
						(rtd\sfx) at (.25,1.75) {};
						
						\node[rt actual] (rtao\sfx) at (1.85,3.5) {};
						\node[rt actual] (rtbo\sfx) at (4,1.75) {};
						\node[rt actual] (rtco\sfx) at (.75,3.5) {};
						\node[rt actual] (rtdo\sfx) at (-.7,1.75) {};
						\draw[rt represent] (rta\sfx) -- (rtao\sfx);
						\draw[rt represent] (rtb\sfx) -- (rtbo\sfx);
						\draw[rt represent] (rtc\sfx) -- (rtco\sfx);
						\draw[rt represent] (rtd\sfx) -- (rtdo\sfx);
						\draw[rt portion] (rtao\sfx)
						.. controls (2.3,4.45) and (4.95,2.4) .. (rtbo\sfx);
						\draw[rt portion] (rtco\sfx)
						.. controls (.3,4.3) and (-1.55,2.3) .. (rtdo\sfx);
						\node[rt note] at (3.85,3.55) {$\Gamma$};
						\node[rt note] at (-.95,3.3) {$R$};
					\end{scope}
				}
				
				\node at (1.55,.05) {(a) Before merging};
				\node at (7.95,.05) {(b) After merging};
				
				\draw[rt formal] (rtaA) -- (rtbA);
				\draw[rt formal] (rtcA) -- (rtdA);
				
				\draw[rt formal] (rtaB) -- (rtcB);
				\draw[rt formal] (rtbB) -- (rtdB);
				
				\draw[rt formal] (-1.2,-.8) -- (-.35,-.8);
				\node[rt note,anchor=west] at (-.2,-.8) {auxiliary edge};
				\draw[rt represent] (3.2,-.8) -- (4.05,-.8);
				\node[rt note,anchor=west] at (4.2,-.8) {edge of $G$};
				\draw[rt portion] (7.1,-.8) -- (7.95,-.8);
				\node[rt note,anchor=west] at (8.1,-.8) {unchanged portion};
			\end{tikzpicture}
			\caption{Merging an alternating cycle $R=R_{j,i}$ into $\Gamma$; the
				shaded box is $W_j$, the dashed line separates $A_j$ from $B_j$, and
				only the auxiliary edges inside $W_j$ are drawn. The auxiliary edges
				$uv$ of $\Gamma$ and $zz'$ of $R$ are replaced by $uz$ and $vz'$. The
				edges of $G$ and the grey portions are unchanged.}
			\label{fig:splice}
		\end{figure}
		
		After all the cycles $R_{j,i}$ have been merged, $\Gamma$ satisfies
		(a), and its vertex set is the disjoint union of $V(\Gamma_0)$ and the
		sets $V(R_{j,i})$. Hence $\Gamma$ has $m_0+h$ auxiliary edges, and (b)
		follows from \eqref{eq:budget}, as
		$m_0+h\le\lambda n/4<\beta\eta n\le\beta|V|$ for every cluster $V$. For
		(c), let $j\in[\ell]$, and for an alternating cycle $\Gamma'$, let
		$a^j_{\Gamma'}$ and $b^j_{\Gamma'}$ be the numbers of auxiliary edges of
		$\Gamma'$ inside $A_j$ and inside $B_j$, respectively. Every vertex of
		$\Gamma'$ lies in exactly one auxiliary edge of $\Gamma'$, and every
		auxiliary edge of $\Gamma'$ that meets $W_j$ lies inside $W_j$, so
		$a^j_{\Gamma'}-b^j_{\Gamma'}=\bigl(|V(\Gamma')\cap A_j|-|V(\Gamma')\cap
		B_j|\bigr)/2$, which is additive over the disjoint union above. Hence,
		by \eqref{eq:Rcount},
		\[
		a^j_{\Gamma}-b^j_{\Gamma}=p_{A_j}-p_{B_j}+\sum_{i=1}^{d_j}1
		=p_{A_j}-p_{B_j}+d_j=|A_j|-|B_j|.
		\]
		Thus $G$ is Hamiltonian by Claim~\ref{cl:skeleton}.
	\end{proof}

\section{Proof of Proposition~\ref{prop:petersen}}\label{sec:example}

Chv\'atal's construction~\cite[Section~5]{Chvatal1973} replaces the
vertices of the Petersen graph by triangles; we replace the vertices of
the Petersen graph $R$ by almost complete graphs instead. We use the following
properties of $R$: it is $3$-regular, has $10$ vertices and $15$
edges, is $3$-connected, and is not Hamiltonian.

\begin{proof}
Let $m\ge25$ be an odd integer, and let $V(R)=[10]$. For each
$i\in[10]$, let $R_i$ be obtained from a copy of $K_m$ with three
distinct vertices $a_i,b_i,c_i$ by deleting the edges $a_ib_i$,
$b_ic_i$, $c_ia_i$ and the edges of a perfect matching on  the remaining
$m-3$ vertices, which exists since $m-3$ is even. We take
$R_1,\ldots,R_{10}$ pairwise vertex-disjoint, and we call $a_i,b_i,c_i$
the \emph{attachment vertices} of $R_i$. Since $R$ is $3$-regular, we
may assign, for each $i$, the three attachment vertices of $R_i$
bijectively to the three edges of $R$ incident with $i$. Let $G$ be
obtained from $R_1\cup\cdots\cup R_{10}$ by adding, for every edge $ij$
of $R$, the edge joining the attachment vertices of $R_i$ and of $R_j$
assigned to $ij$. We call these $15$ edges the \emph{cross edges}.
Every attachment vertex lies in exactly one cross edge, so the cross
edges form a matching.

A vertex of $R_i$ that is not an attachment vertex has degree $m-2$ in
$R_i$ and lies in no cross edge. An attachment vertex has degree $m-3$
in $R_i$ and lies in exactly one cross edge. Thus $G$ is
$(m-2)$-regular on $n=10m$ vertices. In $R_i$, every vertex is
nonadjacent to at most two other vertices. Consequently,
\begin{equation}\label{eq:petersen_pieces}
\text{if $X\subseteq V(R_i)$ and $|X|\ge5$, then $R_i[X]$ is connected,}
\end{equation}
since two nonadjacent vertices of $X$ have at least $|X|-4\ge1$ common
neighbors in $X$. In particular, each $R_i$ is connected, and since $R$
is connected, so is $G$. Also,
\[
\frac{\delta(G)}{|V(G)|}=\frac{m-2}{10m}\ge\alpha
\quad\Longleftrightarrow\quad
m\ge\frac{2}{1-10\alpha},
\]
and there are arbitrarily large odd integers $m$ satisfying this
inequality.

Suppose that $C$ is a Hamiltonian cycle of $G$. For each $i$, the cycle
$C$ has vertices both in $V(R_i)$ and outside it, so $C$ contains a
positive even number of edges between $V(R_i)$ and its complement.
These edges are cross edges, and only three cross edges meet $V(R_i)$.
Hence $C$ contains exactly two of them. This implies that $C\cap R_i$
is a Hamiltonian path of $R_i$. Therefore, contracting each $R_i$
into a single vertex 
turns $C$ into a Hamiltonian cycle of $R$, a contradiction. Thus $G$ is
not Hamiltonian.

It remains to show that $G$ is $3/2$-tough. Let $S\subseteq V(G)$
with $c(G-S)\ge2$, and let $B$ be the set of indices $i$ with
$|V(R_i)\setminus S|\le4$.

Assume first that $B\ne\emptyset$. For $i\notin B$, the graph
$R_i-S$ is connected by \eqref{eq:petersen_pieces}, so at most
$10-|B|$ components of $G-S$ contain a vertex of such a graph. Every
other component of $G-S$ is contained in the set
$\bigcup_{i\in B}(V(R_i)\setminus S)$, which has at most $4|B|$
vertices. Hence
\[
c(G-S)\le(10-|B|)+4|B|=10+3|B|.
\]
On the other hand, $|S|\ge(m-4)|B|\ge21|B|$. Since $|B|\ge1$, we have
$21|B|\ge\frac32(10+3|B|)$, and so $|S|\ge\frac32c(G-S)$.

Assume now that $B=\emptyset$; we follow the arguments in the proofs
of~\cite[Theorems~4.1 and~5.1]{Chvatal1973}. Every $R_i-S$ is then nonempty and,
by \eqref{eq:petersen_pieces}, connected. Let $F$ be the set of edges
of $R$ whose cross edges have an end in $S$. Two graphs $R_i-S$ and
$R_j-S$ are joined by an edge of $G-S$ if and only if $ij\in
E(R)\setminus F$, so the components of $G-S$ correspond to the
components of $R-F$. Let $U_1,\ldots,U_s$ be the vertex sets of the
components of $R-F$, so that $s=c(G-S)\ge2$. Since the cross edges form
a matching, distinct edges of $F$ give distinct vertices of $S$, and so
$|S|\ge|F|$. Since $s\ge2$, each $U_h$ is a nonempty proper subset of
$V(R)$, so, as $R$ is $3$-connected, at least three edges of $R$ join
$U_h$ to its complement.
These edges join distinct components of $R-F$ and therefore lie in
$F$, and each edge of $F$ is counted in this way for at most two sets
$U_h$. Hence $2|F|\ge3s$, and so $|S|\ge|F|\ge\frac32c(G-S)$. Thus $G$
is $3/2$-tough.
\end{proof}

\section*{Statement on the use of AI}

The research was directed by the author, who posed the problem and chose
the method, namely the robust partition theorem of K\"uhn, Lo, Osthus,
and Staden~\cite{KLOS2015}. ChatGPT Pro, through its chat interface,
was used as a research assistant in many rounds of discussion with the
author. Over these
discussions, it produced a first version of the proof, for $6$-tough
graphs, and, at the author's direction, the exchange argument
in the proof of Lemma~\ref{lem:twotrees}, which replaced a topological
step of the first version and lowered the toughness requirement to $4$. The proof presented here was written by the author,
who rewrote the whole argument and simplified most of its parts. In
particular, the author generalized Lemma~\ref{lem:twotrees} from two
edge-disjoint spanning trees in $4$-tough graphs to $k$ edge-disjoint
spanning trees under a condition implied by $2k$-toughness, gave the
short proof of Lemma~\ref{lem:unbalanced}, and introduced the alternating
cycles and the cycle-merging argument in the proof of
Theorem~\ref{thm:main}. Claude Opus 5.5 was used, under the author's
direction, to implement revisions, proofread the draft, and check the
references. The author has verified all arguments and takes full
responsibility for the content of this paper.

\end{document}